\documentclass[11pt]{article}
\usepackage{amsmath,amssymb,amsthm}
\usepackage[margin=1in]{geometry}
\usepackage{graphicx}
\usepackage{color}
\usepackage{caption,subcaption}
\usepackage{mathtools}
\usepackage{authblk}
\usepackage[sort&compress,numbers]{natbib}  
\usepackage{hyperref}

\newcommand{\RR}{\mathbb{R}}
\newcommand{\CC}{\mathbb{C}}

\newcommand{\OO}[1]{\mathcal{O}\paren{#1}}
\newcommand{\paren}[1]{\left( #1 \right)}

\newcommand{\wh}[1]{\widehat{#1}}
\newcommand{\wt}[1]{\widetilde{#1}}

\newcommand{\abs}[1]{\left\lvert#1\right\rvert}

\newcommand{\Schrod}{Schr\"{o}dinger\ }
\newcommand{\legt}{\mathcal{T}}
\newcommand{\tmax}{t_{\max}}

\renewcommand{\Re}{\operatorname{Re}}
\renewcommand{\Im}{\operatorname{Im}}

\DeclareMathOperator\erfc{erfc}
\DeclareMathOperator\erfi{erfi}

\title{A fast, high-order numerical method for the simulation of
single-excitation states in quantum optics}

\author[1]{Jeremy Hoskins}
\author[2,3]{Jason Kaye}
\author[2]{Manas Rachh}
\author[4]{John C. Schotland}

\affil[1]{{\footnotesize Department of Statistics, University of Chicago, Chicago,
IL 60637, USA}}
\affil[2]{{\footnotesize Center for Computational Mathematics, Flatiron Institute, New York, NY 10010, USA}}
\affil[3]{{\footnotesize Center for Computational Quantum Physics, Flatiron Institute, New York, NY 10010, USA}}
\affil[4]{{\footnotesize Department of Mathematics and Department of Physics, Yale University, New Haven, CT
06511, USA}}

\date{}

\begin{document}

\maketitle

\begin{abstract}

We consider the numerical solution of a nonlocal partial differential
equation which describes the phenomenon of collective spontaneous emission in
a two-level atomic system containing a single photon. We reformulate the
problem as an integro-differential equation for the atomic degrees of
freedom, and describe an efficient solver for the case of a Gaussian
atomic density. The problem of history dependence arising from the
integral formulation is addressed using sum-of-exponentials history
compression. We demonstrate the solver on two systems of physical
interest: in the first, an initially-excited atom decays into a photon
by spontaneous emission, and in the second, a photon pulse is used to an
excite an atom, which then decays.

\end{abstract}

\textbf{\textit{Keywords ---}} {\small quantum optics; nonlocal partial
differential equations; Volterra integro-differential equations;
sum of exponentials compression}

\section{Introduction}
Many-body problems in quantum optics are of interest in the study of
cold-atom systems, quantum waveguides, and quantum semiconductor devices, among others, with applications to quantum computing, quantum information processing, 
and precision measurements~\cite{Haroche_2006,Gardiner_2015,Liao_2016,Roy_2017,Kira_2011,Kimble_2008,Riedmatten_2008,Bloch_2012}.
The simplest such problem arises in a system of
two-level atoms interacting with a single photon. In this setting, the
propagation of a single-photon state is governed by the system of
partial differential equations~\cite{kraisler22}
\begin{equation}
\begin{aligned}
  \label{eq:system}
    i \partial_t u(x,t) &= c (-\Delta)^{1/2} u(x,t) +
    g \rho(x) a(x,t), \quad (x,t)\in {\mathbb R}^{d+1}, \\ 
    i \partial_t a(x,t) &= \Omega a(x,t) + g u(x,t).
\end{aligned}
\end{equation}
Here $u$ is the probability amplitude for creating a photon, $a$ is the
probability amplitude for exciting an atom, $\rho$ is the atomic number density, $\Omega$ is the atomic resonance frequency, and $g$ is the
atom-field coupling constant. The amplitudes obey the normalization
condition
\begin{equation} \label{eq:normalization}
  \int_{\RR^d} \paren{|u(x,t)|^2 + \rho(x)|a(x,t)|^2} \, dx = 1,
\end{equation}
which has the interpretation that $|u|^2$
is the one-photon probability density and that $\rho|a|^2$ is the atomic
probability density. In physical terms, (\ref{eq:system}) describes the
process of collective spontaneous emission. That is, suppose that an
atom is initially in its excited state and there are no photons present
in the field. The atom can then decay, transferring its excitation to
the field, which can then excite the remaining atoms, causing them to
decay in a similar manner and so on.

Eq.~\eqref{eq:system} has been investigated in several cases of
interest, including a single atom, a uniform medium of constant density,
and a statistically homogeneous random medium~\cite{kraisler22}. This
paper is the first in a series devoted to the analysis and numerical solution of
\eqref{eq:system}.
We note that standard numerical methods are not readily applicable to
this problem, which was originally introduced in Ref. \cite{kraisler22},
and to our knowledge this is the first paper which discusses its
numerical solution.
In order to illustrate the difficulty, we outline the drawbacks of two
possible approaches. 

\paragraph{Physical domain discretization}
We could consider discretizing the first equation in
\eqref{eq:system} directly in physical space using a finite
difference or finite element method, and then solve the resulting system of ODEs. However, the nonlocal character of
the fractional Laplacian operator $(-\Delta)^{1/2}$, which is given by
\[(-\Delta)^{1/2} f(x) =
\frac{\Gamma\paren{\frac{d+1}{2}}}{\pi^{\frac{d+1}{2}}}
\int_{\mathbb{R}^d} \frac{f(x)-f(y)}{\abs{x-y}^{d+1}} \, dy, \]
leads to two related difficulties. First, any discretization of the
operator would produce a dense matrix, leading to a large cost per time step
in the absence of suitable fast algorithms. Perhaps more importantly,
the photon field $u(x,t)$ would need to be discretized on a domain
containing its
full numerical support, which spreads rapidly. This would, in practice,
limit simulations to very short  times. One possible remedy 
would be to truncate the computational domain and impose suitable
artificial outgoing boundary conditions, but for large systems the cost of discretizing the photon field in the truncated computational domain would remain an issue.

\paragraph{Fourier domain discretization}
The above observations suggest working in the Fourier domain, in which the
action of the fractional Laplacian is diagonal: 
\[(-\Delta)^{1/2} f(x) = \frac{1}{(2\pi)^d} \int_{\RR^d} e^{i
\xi \cdot x} \abs{\xi} \wh{f}(\xi) \, d \xi,\]
where $\wh{f}(\xi)$ is the Fourier transform of $f$, which is defined by
\[\wh{f}(\xi) = \int_{\RR^d} e^{-i \xi \cdot x} f(x) \, dx.\]
One could design a Fourier pseudospectral method, such that at each time
step, the action of the fractional Laplacian is computed in the Fourier
domain, and the product $\rho(x) a(x,t)$ is computed in the physical
domain. Such methods are commonly used to solve PDEs of evolution, such
as the time-dependent \Schrod
equation, involving a Laplacian term diagonal in the Fourier domain, and
a second term which is more easily computed in the physical domain
\cite{boyd01,blanes00}. Here, we encounter the Fourier domain
manifestation of the same problem. Namely, spreading in the physical
domain corresponds to oscillation in the Fourier domain, and we
obtain a photon amplitude which becomes more and more oscillatory in the
Fourier domain as time progresses. As a result, one would expect 
the computational cost to scale at least quadratically with the
propagation time.

\vspace{\baselineskip}

Our approach is to recast \eqref{eq:pde} as a Volterra integral equation 
for the atomic amplitude. In particular, we eliminate the photon field using a suitable Green's function, obviating
the need to discretize large spatial domains. The number of degrees of
freedom in the required discretization 
depends only on the size of the support of $\rho$. As such, our method enables fast and accurate 
simulations over long times.

We begin by constructing the Green's function
for the homogeneous part of the equation describing $u$, which satisfies
\begin{equation}
\begin{aligned}
  \label{eq:gpde}
    i \partial_t G(x,t) &= c (-\Delta)^{1/2} G(x,t) \\
    \lim_{t \to 0^+} G(x,t) &= \delta(x).
\end{aligned}
\end{equation}
The solution in the Fourier domain is given by
\begin{equation} \label{eq:gfun}
  \wh{G}(\xi,t) = e^{-ic\abs{\xi}t}.
\end{equation}
This implies that in the case $g = 0$, $u(\xi,t)$ is given by
\[u(\xi,t) = \wh{G}(\xi,t) \wh{u_0}(\xi) = e^{-ic\abs{\xi}t}
\wh{u_0}(\xi),\]
from which the oscillatory behavior is clear.

We wish to make use of the Green's function representation
of $u(x,t)$, but to avoid discretizing it in the Fourier domain. To proceed,
we rewrite \eqref{eq:system} as
\begin{equation}
\begin{aligned}
  \label{eq:pde}
    i \partial_t u(x,t) &= c (-\Delta)^{1/2} u(x,t) +
    \frac{g}{\sigma^d} \rho(x/\sigma) a(x,t) ,  \\ 
    i \partial_t a(x,t) &= \Omega a(x,t) + g u(x,t) , \\
    u(x,0) &= u_0(x) , \\
    a(x,0) &= a_0(x) ,
\end{aligned}
\end{equation}
where the density $\rho$ has been rescaled by the length $\sigma$, which
characterizes the spatial extent of the atoms. Next we reformulate
\eqref{eq:pde} as a Volterra
integro-differential equation in the unknown $b(x,t) =
\frac{\rho(x/\sigma)}{\sigma^d} a(x,t)$ alone. Applying the Duhamel principle to the first equation in \eqref{eq:pde} gives
\begin{equation} \label{eq:uint}
  u(x,t) = \int_{\RR^d} G(x-y,t) u_0(y) \, dy - i \frac{g}{\sigma^d} \int_0^t
\int_{\RR^d} G(x-y,t-s) \rho(y/\sigma) a(y,s) \, dy \, ds.
\end{equation}
Substituting the above into the
second equation in \eqref{eq:pde} and multiplying by
${\rho(x/\sigma)}/{\sigma^d}$
gives
\begin{equation} \label{eq:vide}
  \partial_t b(x,t) = -i \Omega b(x,t) - g^2 \frac{\rho(x/\sigma)}{\sigma^d} \int_0^t \int_{\RR^d} G(x-y,t-s)
  b(y,s) \, dy \, ds - ig \frac{\rho(x/\sigma)}{\sigma^d} U(x,t),
\end{equation}
where we have defined 
\[U(x,t) = \int_{\RR^d} G(x-y,t) u_0(y) \, dy,\]
which is the free evolution of the photon amplitude $u_0(x)$.
If \eqref{eq:vide} is solved, the photon amplitude can be recovered
as a matter of post-processing using \eqref{eq:uint}.

The main
advantage of solving \eqref{eq:vide} over the formulations mentioned
above is that for a localized density $\rho(x)$, $b(x,t)$ remains
localized as well. The price we pay is a dense dependence of the
solution $b(x,t)$ on its history $b(x,s)$ for $0 \leq s < t$.  Indeed, it
appears that each time step, we must evaluate the history
integral on the right hand side of \eqref{eq:vide}. This leads to an
algorithm which, for a given accuracy, has a computational cost scaling
as $\OO{N^2}$ in the number $N$ of time steps, and a memory requirement scaling
as $\OO{N}$. This is a typical challenge associated with the application of Volterra
integral operators, and several techniques have been proposed to
address it, particularly in the context of solving Volterra integral equations
\cite{hairer85,veerapaneni07,jiang15,wang19,kaye21_2,kaye22_tdse,dolz21} and
applying Volterra integral operators corresponding to nonlocal
transparent boundary conditions
\cite{alpert00,lubich02,jiang04,schadle06,jiang08,kaye18}. We will make use of one such
approach -- the sum of exponentials approximation method -- to obtain a
high-order accurate numerical method with $\OO{N \log N}$ computational
complexity and $\OO{\log N}$ memory complexity.

We focus in this article on the case of a Gaussian atomic density in one
spatial dimension. There is no fundamental difficulty in extending
our method to densities comprised of sums of Gaussians, and to three
spatial dimensions. These extensions will be addressed in a forthcoming
publication. A generalization to other densities may also be possible,
but Gaussian and sum-of-Gaussian densities are a suitable physical model
for many systems of contemporary interest. We will see
that the present case already exhibits nontrivial dynamics which
are expected to appear in three dimensions as well.

This article is organized as follows. In Section \ref{sec:setup}, we
describe the mathematical setup for our numerical method. We describe our
high-order time-stepping algorithm in Section \ref{sec:nummethod}, and
fill in technical details involving the representation and evaluation of
certain special functions in Section \ref{sec:kernels}. In Section \ref{sec:results}
we present numerical results which demonstrate the accuracy of the method
and give insight into the behavior of the solution for two physically
meaningful examples. Section \ref{sec:conclusion} concludes with a
discussion of several open questions and future research directions.

\section{Problem setup} \label{sec:setup}

To set up our numerical method we will
represent the atom amplitude $a(x,t)$ in the one-dimensional case by an expansion
\begin{equation} \label{eq:modalexp}
  a(x,t) = \sum_{n=0}^{p-1} a_n(t) f_n(x/\sigma).
\end{equation}
Here $\{f_n(x)\}_{n=0}^{p-1}$ are the first $p$ polynomials orthonormal with
respect to $\rho(x)$, so that $\{f_n(x/\sigma)\}_{n=0}^{p-1}$ are
orthonormal with respect to the scaled density ${\rho(x/\sigma)}/{\sigma}$.
We will first derive a coupled set
of Volterra integral equations (VIEs) for the modal coefficients
$a_n(t)$. We will then obtain explicit
expressions for the case in which the atomic density $\rho$ is a
Gaussian.
Finally, we will show how to recover the photon amplitude from the
coefficients $a_n$ of the atom amplitude.

\subsection{Volterra integral equation for the atomic degrees of
freedom}

Substituting \eqref{eq:modalexp} into \eqref{eq:vide}, integrating against
$f_m(x/\sigma)$, and defining 
\[U_m(t) = \frac{1}{\sigma} \int_{-\infty}^\infty \rho(x/\sigma)
f_m(x/\sigma) U(x,t) \,
dx,\]
we obtain
\begin{multline*}
\dot{a}_m(t) = -i \Omega a_m(t) \\
  - \frac{g^2}{\sigma^2} \sum_{n=0}^{p-1} \int_{-\infty}^\infty
  \rho(x/\sigma) f_m(x/\sigma) \int_0^t a_n(s) \int_{-\infty}^\infty G(x-y,t-s)
  \rho(y/\sigma) f_n(y/\sigma) \, dy \, ds \, dx - i g U_m(t),
\end{multline*}
where the dot denotes a derivative with respect to time.
From \eqref{eq:gfun}, we have
\[G(x,t) = \frac{1}{2 \pi} \int_{-\infty}^\infty e^{i (\xi x - c
\abs{\xi} t)} \, d\xi,\]
which gives
\begin{align*}
  &\int_{-\infty}^\infty \rho(x/\sigma) f_m(x/\sigma) \int_0^t a_n(s) 
\int_{-\infty}^\infty 
  G(x-y,t-s) \rho(y/\sigma) a_n(s) f_n(y/\sigma) \, dy \, ds \, dx \\
  &= \frac{1}{2\pi} \int_0^t a_n(s) \int_{-\infty}^\infty e^{-ic\abs{\xi}(t-s)}
  \paren{\int_{-\infty}^\infty e^{i \xi x} \rho(x/\sigma) f_m(x/\sigma)
  \, dx}
  \paren{\int_{-\infty}^\infty e^{-i \xi y} \rho(y/\sigma) f_n(y/\sigma)
  \, dy} \, d\xi \, ds\\
  &= \frac{\sigma^2}{2\pi} \int_0^t a_n(s) \int_{-\infty}^\infty
  e^{-ic\abs{\xi}(t-s)} \wh{\paren{\rho f_m}}(-\sigma \xi) \wh{\paren{\rho
  f_n}}(\sigma \xi) \, d\xi \, ds \\
  &= \frac{\sigma^2}{2\pi} \int_0^t a_n(s) \int_{0}^\infty
  e^{-i c \xi (t-s)} \Phi_{mn}(\sigma \xi) \, d\xi \, ds ,
\end{align*}
where 
\[\Phi_{mn}(\xi) = \phi_m(\xi) \phi_n(-\xi) + \phi_m(-\xi) \phi_n(\xi)\]
with
\[\phi_n(\xi) = \wh{\paren{\rho f_n}}(\xi).\]
Defining
\[J_{mn}(t) = \int_0^\infty e^{-i \xi t} \Phi_{mn}(\xi) \, d\xi ,\]
we obtain
\[\dot{a}_m(t) = -i \Omega a_m(t) - \frac{g^2}{2\pi\sigma}
\sum_{n=0}^{p-1}
  \int_0^t J_{mn}\paren{\frac{c}{\sigma}(t-s)} a_n(s) \, ds - i g U_m(t).\]
The change of variables 
\begin{equation} \label{eq:cov}
  \alpha_m(t) = e^{i \Omega t} a_m(t)
\end{equation}
gives
\[\dot{\alpha}_m(t) = - \frac{g^2}{2\pi\sigma} \sum_{n=0}^{p-1} \int_0^t
e^{i \Omega (t-s)} J_{mn}\paren{\frac{c}{\sigma}(t-s)} \alpha_n(s) \, ds - i g e^{i \Omega t} U_m(t).\]
Integrating both sides in time and swapping the
order of integration yields
\begin{equation} \label{eq:vie}
  \alpha_m(t) + \frac{g^2}{2\pi c} \sum_{n=0}^{p-1} \int_0^t
  K_{mn}\paren{\frac{c}{\sigma} (t-s)} \alpha_n(s) \, ds = a_m(0) - i g \int_0^t e^{i \Omega s}
U_m(s) \, ds
\end{equation}
with
\begin{equation} \label{eq:kappadef}
  K_{mn}(t) = \int_0^t e^{i \frac{\Omega \sigma}{c} s} J_{mn}(s) \, ds.
\end{equation}
The above is a collection of coupled second-kind VIEs for $\alpha_m(t)$, $m =
0,\ldots,p-1$, from which $a(x,t)$ can be
recovered using \eqref{eq:modalexp} and \eqref{eq:cov}.

We pause to consider the calculation of the total probability, given as
in \eqref{eq:normalization} by
\begin{equation} \label{eq:totprob}
1 = \frac{1}{\sigma} \int_{-\infty}^\infty
\abs{a(x,t)}^2 \rho(x/\sigma) \, dx + \int_{-\infty}^\infty
  \abs{u(x,t)}^2 \, dx \equiv P_a(t) + P_u(t).
\end{equation}
Here, we have defined $P_a$ and $P_u$ as the atomic and
photonic contributions to the
total probability, respectively.
It is straightforward to calculate $P_a$, a quantity of physical
interest, within our framework: 
\begin{equation} \label{eq:energya}
  \begin{aligned}
    P_a(t) &= \frac{1}{\sigma} \int_{-\infty}^\infty \abs{a(x,t)}^2
    \rho(x/\sigma) \, dx \\
    &= \frac{1}{\sigma} \sum_{n=0}^{p-1} \sum_{m=0}^{p-1} 
    a_m^*(t) a_n(t) \int_{-\infty}^\infty f_m^*(x/\sigma)
    f_n^*(x/\sigma) \rho(x/\sigma) \, dx \\
    &= \sum_{n=0}^{p-1} \abs{a_n(t)}^2 = \sum_{n=0}^{p-1}
    \abs{\alpha_n(t)}^2.
  \end{aligned}
\end{equation}

\subsection{Gaussian atomic density}

Let us take the atomic density to be a Gaussian,
\[\rho(x) = \frac{e^{-x^2}}{\sqrt{\pi}}.\]
Then
\[f_n(x) = \frac{H_n(x)}{\sqrt{2^n n!}}, \] 
with $H_n$ the Hermite polynomial of degree $n$, defined by
\begin{equation} \label{eq:hermitedef}
  H_n(x) e^{-x^2} = (-1)^n \frac{d^n}{dx^n} e^{-x^2}.
\end{equation}
The above follows from the formula \cite[Eqn. 7.374.1]{gradshteyn07}
\[\int_{-\infty}^\infty H_m(x) H_n(x) e^{-x^2} \, dx = \sqrt{\pi} 2^n n!
\delta_{mn}.\]
Taking the Fourier transform of
\eqref{eq:hermitedef} gives
\begin{equation} \label{eq:phindef}
\phi_n(\xi) = \wh{\paren{\rho f_n}}(\xi) = \frac{(-i)^n}{\sqrt{2^n n!}} \xi^n e^{-\xi^2 / 4}.
\end{equation}
In particular, we find that $\Phi_{mn}(\xi) = 0$ if $m$ is even and $n$ is
odd or vice versa, and otherwise
\[\Phi_{mn}(\xi) = 2 (-1)^m \phi_m(\xi) \phi_n(\xi) = \frac{(-1)^m
(-i)^{m+n}}{\sqrt{2^{m+n-2} m! n!}} \xi^{m+n}
e^{-\xi^2 / 2}.\]
We remark that the vanishing of $\Phi_{mn}$ for odd $m+n$ is a consequence of the 
symmetry of $\rho$. For more general densities, all $\Phi_{mn}$ will be non-zero.
The kernel $J_{mn}$ is then given by
\[J_{mn}(t) = \frac{(-1)^m (-i)^{m+n}}{\sqrt{2^{m+n-2} m! n!}}
 \int_0^\infty \xi^{m+n} e^{-\xi^2 / 2 - i \xi t} \, d\xi\]
if $m$ and $n$ are even or odd together, and zero otherwise.
We define
\begin{equation} \label{eq:jndef}
  j_n(t) = \frac{2}{\Gamma\paren{\frac{n+1}{2}}} \int_0^\infty \xi^n
e^{-\xi^2 - i \xi t} \, d\xi.
\end{equation}
Here $\Gamma$ is the Gamma function, and the normalization is chosen so
that $j_n(0) = 1$. A change of variables gives
\[J_{mn}(t) = 
\begin{cases}
(-1)^m (-i)^{m+n}
\frac{\Gamma\paren{\frac{m+n+1}{2}}}{\sqrt{\frac{m! n!}{2}}}
  j_{m+n}\paren{\sqrt{2} t} & \text{if } m+n \equiv 0 \pmod{2} \\
  0 & \text{otherwise.}
\end{cases}
\]
We also define
\begin{equation} \label{eq:kndef}
k_n(t) = \int_0^t e^{i \frac{\Omega \sigma}{c} s} j_n(\sqrt{2} s) \, ds
\end{equation}
so that
\begin{equation} \label{eq:kmndef}
K_{mn}(t) = 
\begin{cases}
(-1)^m (-i)^{m+n}
\frac{\Gamma\paren{\frac{m+n+1}{2}}}{\sqrt{\frac{m! n!}{2}}}
  k_{m+n}(t) &\text{if } m+n \equiv 0 \pmod{2} \\
  0 &\text{otherwise.}
\end{cases}
\end{equation}

\subsection{Recovering the photon amplitude}

The photon amplitude is given by \eqref{eq:uint}. The first term,
$U(x,t)$,
describes the contribution to the amplitude of the initial photon field
configuration, and is straightforward to
compute by Fourier transform as long as $u_0$ is well-behaved.

For the second term, we write
\begin{align*}
  u(x,t) - U(x,t) &= -\frac{ig}{\sigma} \int_0^t \int_{-\infty}^\infty 
  G(x-y,t-s) \rho(y/\sigma) a(y,s) \, dy \, ds \\
  &= -\frac{ig}{2 \pi \sigma} \int_0^t \int_{-\infty}^\infty
  e^{-ic\abs{\xi}(t-s)} e^{i \xi x} \int_{-\infty}^\infty
  e^{-i \xi y} \rho(y/\sigma) a(y,s) \, dy \, d\xi \, ds \\
  &= -\frac{ig}{2 \pi \sigma} \sum_{n=0}^{p-1} \int_0^t a_n(s) \int_{-\infty}^\infty 
  e^{-ic\abs{\xi}(t-s)} e^{i \xi x} \int_{-\infty}^\infty 
  e^{-i \xi y} \rho(y/\sigma) f_n(y/\sigma) \, dy \, d\xi \, ds \\
  &= -\frac{ig}{2 \pi \sigma} \sum_{n=0}^{p-1} \int_0^t a_n(s) \int_{-\infty}^\infty 
  e^{-ic\abs{\xi}(t-s)/\sigma} e^{i \xi x/\sigma} \phi_n(\xi) \, d\xi \,
  ds\\
  &= -\frac{ig}{2 \pi \sigma} \sum_{n=0}^{p-1} \int_0^t a_n(s) \int_0^\infty
  e^{-ic \xi (t-s)/\sigma} \paren{e^{i \xi x/\sigma} \phi_n(\xi) + e^{-i
  \xi x/\sigma} \phi_n(-\xi)} \, d\xi \, ds.
\end{align*}
Once we have solved \eqref{eq:vie}, we can recover the coefficients
$a_m(t)$ from \eqref{eq:cov}, and compute the photon amplitude as
above. In the case of a Gaussian atomic density,
\eqref{eq:phindef} and \eqref{eq:jndef} yield
\begin{align*}
  u(x,t) &- U(x,t) \\
  &= -\frac{ig}{2 \pi \sigma} \sum_{n=0}^{p-1}
  \frac{(-i)^n}{\sqrt{2^n n!}} \int_0^t a_n(s) \int_0^\infty
  e^{-ic \xi (t-s)/\sigma-\xi^2/4} \paren{e^{i \xi x/\sigma} \xi^n +
  (-1)^n e^{-i \xi x/\sigma} \xi^n} \, d\xi \, ds \\
  &= \begin{multlined}[t][.83\textwidth]
    -\frac{ig}{2 \pi \sigma} \sum_{n=0}^{p-1} 
  \frac{(-i)^n 2^\frac{n}{2} \Gamma\paren{\frac{n+1}{2}}}{\sqrt{n!}}  \int_0^t
    a_n(s) \left[j_n\paren{\frac{2\paren{c(t-s)-x}}{\sigma}} \right. \\
    \left.+ (-1)^n j_n\paren{\frac{2\paren{c(t-s)+x}}{\sigma}}\right] \, ds , 
  \end{multlined}
\end{align*}
after some manipulation.

\section{Discretization and numerical solution} \label{sec:nummethod}

We use a high-order implicit Gauss-Legendre collocation
method to discretize and solve the VIE \eqref{eq:vie}. As is typical with VIEs, the primary computational
bottleneck is the evaluation of history integrals at each time step. The
naive cost of these evaluations scales quadratically with the
total number of time steps, but we will show that it can reduced by
splitting the history integrals into local and history parts, and
deriving recurrences for the latter using sum-of-exponentials representations of the
kernels $K_{mn}(t)$.

We begin by describing our discretization scheme. We divide the time
interval $[0,T]$ into $N$ uniform subintervals $\{[(j-1) \Delta t,j \Delta
t]\}_{j=0}^{N-1}$, with $\Delta t = T/N$. Let $\{\tau_k\}_{k=0}^{q-1}$ be
the collection of $q$ Gauss-Legendre nodes, rescaled and shifted to the interval
$[0,\Delta t]$. We place $q$ Gauss-Legendre nodes on each subinterval,
so that the full set of collocation nodes is given by $t_{jk} = (j-1) \Delta
t + \tau_k$ for $j=1,\ldots,N$ and $k=0,\ldots,q-1$.

We denote the numerical
approximation of $\alpha_m(t_{jk})$ by $\alpha_{m,j,k}$. In
addition to this so-called grid representation of the numerical solution, we
will also sometimes represent the numerical solution on a subinterval
$[(j-1) \Delta t, j \Delta t]$ by
\begin{equation} \label{eq:legrep}
  \alpha(t) \approx \sum_{k=0}^{q-1} \wh{\alpha}_{m,j,k} P_k^j(t),
\end{equation}
where $P_k^j(t)$ is the Legendre polynomial of degree $k$ on the interval $[(j-1) \Delta t, j \Delta t]$; that is, $P_k^j(t) = P_k(t - (j-1)\Delta t)$, where $P_k(\tau)$ is the Legendre polynomial of degree $k$ on $[0,\Delta t]$. One can
transform back and forth between the grid representation $\alpha_{m,j,k}$ and the Legendre
coefficient representation $\wh{\alpha}_{m,j,k}$ on the $j$th
subinterval by interpolation of the
expansion \eqref{eq:legrep} at the Gauss-Legendre nodes $t_{jk}$.
Indeed, we have
\[\alpha_{m,j,k} = \sum_{l=0}^{q-1} P_l^j(t_{jk}) \wh{\alpha}_{m,j,l} =
\sum_{l=0}^{q-1} P_l(\tau_k) \wh{\alpha}_{m,j,l},\]
and the matrix $\legt_{kl} = P_l(\tau_k)$ is well-conditioned \cite{gautschi83}. We can
therefore obtain the grid representation from the coefficient
representation by applying $\legt$, and the coefficient representation from
the grid representation by applying $\legt^{-1}$. We refer to $\legt^{-1}$ as
the \emph{discrete Legendre transform} matrix.

We split the integral operator in \eqref{eq:vie} into three
pieces:
\begin{align*}
  \int_0^t K_{mn}(t-s) \alpha_n(s) \, ds &= \paren{\int_{(j-1) \Delta
  t}^t + \int_{t_j^*}^{(j-1) \Delta t} + \int_0^{t_j^*}}
  K_{mn}\paren{\frac{c}{\sigma} \paren{t-s}} \alpha_n(s) \, ds \\
  &\equiv C_{m,n,j}(t) + L_{m,n,j}(t) + H_{m,n,j}(t).
\end{align*}
Here, the labels of the three integrals stand for current-time, local, and
history, respectively. We define $t_j^* = \max\paren{0,(j-M)\Delta t}$ for a
fixed positive integer $M \leq N$, which is the number of time steps in
the current and local intervals in the time domain. The
local interval is empty initially, and grows to a maximum length of $(M-1) \Delta t$, whereas
the history interval is empty until $j = M+1$, after which it grows by
$\Delta t$ each time step. The splitting into local and history parts is made because
the sum-of-exponentials representation of $K_{mn}(t)$ is only valid sufficiently far
into the history, and later $M$ will be chosen based on this domain of
validity. The further splitting off of the current time part
is made to conveniently address implicit time-stepping.

To discretize, we use the notation $C_{m,n,j,k} \approx C_{m,n,j}(t_{jk})
\equiv C_{m,n,j}((j-1)\Delta t + \tau_k)$, and similarly for $L_{m,n,j,k}$ and
$H_{m,n,j,k}$. Then 
rearranging and evaluating at $t = t_{jk}$, the
discretization of the VIE \eqref{eq:vie} can be written as
\begin{equation} \label{eq:viedisc}
  \alpha_{m,j,k} + \frac{g^2}{2\pi c} \sum_{n=1}^p C_{m,n,j,k} = -
\frac{g^2}{2\pi c} \sum_{n=1}^p \paren{L_{m,n,j,k} + H_{m,n,j,k}} + f_m(t_{jk})
\end{equation}
where we consider $f_m(t) = a_m(0) - i g \int_0^t e^{i \Omega s} U_m(s)
\, ds$ as a known source term. We note that at a given time step
$j_0$, all of the quantities on the right-hand side depend on the
numerical solution $\alpha_{m,j,k}$ computed only in the first $j_0-1$
time steps, $0 \leq j \leq j_0-1$,
whereas the left hand side depends on the current-time solution
$\alpha_{m,j_0,k}$.

\subsection{The current-time term}

We have
\begin{align*}
  C_{m,n,j,k} &= \int_{(j-1)\Delta t}^{t_{jk}}
  K_{mn}\paren{\frac{c}{\sigma} \paren{t_{jk}-s}}
  \alpha_n(s) \, ds\\
  &= \sum_{l=0}^{q-1} \wh{\alpha}_{n,j,l} \int_{(j-1)\Delta t}^{t_{jk}}
  K_{mn}\paren{\frac{c}{\sigma} \paren{t_{jk}-s}} P_l^j(s) \, ds \\
  &= \sum_{l=0}^{q-1} \wh{\alpha}_{n,j,l} \int_{0}^{\tau_k} 
  K_{mn}\paren{\frac{c}{\sigma} \paren{\tau_k-s}} P_l(s) \, ds \\
  &= \sum_{l=0}^{q-1} \wh{\mathcal{C}}_{m,n,k,l} \wh{\alpha}_{n,j,l},
\end{align*}
where $\wh{\mathcal{C}}_{m,n,k,l} = \int_{0}^{\tau_k}
K_{mn}\paren{\frac{c}{\sigma} \paren{\tau_k-s}}
P_l(s) \, ds$. For each fixed $m$ and $n$,
the array $\wh{\mathcal{C}}_{m,n,k,l}$ can be precomposed with the
discrete Legendre transform matrix $\legt^{-1}$, defined above, to obtain an array $\mathcal{C}_{m,n,k,l}$ with
\begin{equation} \label{eq:cmat}
  C_{m,n,j,k} = \sum_{l=0}^{q-1} \mathcal{C}_{m,n,k,l}
\alpha_{n,j,l}.
\end{equation}
This allows us to work directly with the grid representation of
$\alpha_n$.

The $p^2 q^2$ quantities $\wh{\mathcal{C}}_{m,n,k,l}$, and
therefore $\mathcal{C}_{m,n,k,l}$, can be
precomputed and stored. By plugging the expression \eqref{eq:kmndef} for
$K_{mn}(t)$ in terms of $k_{m+n}(t)$ into
the definition of $\wh{\mathcal{C}}_{m,n,k,l}$, we observe that this can be
accomplished by computing only $p q^2$ integrals -- corresponding to the
different choices of even $m+n$ for $m,n = 0,\ldots,p-1$ -- and
scaling the results by constants depending on $m$ and $n$.

\subsection{The local term} \label{sec:local}

We first split the local term into integrals over the subintervals
defining each time step, and then take a similar approach as for the current-time term:
\begingroup
\allowdisplaybreaks
\begin{align*}
  L_{m,n,j,k} &= \int_{t_j^*}^{(j-1) \Delta t} K_{mn}\paren{\frac{c}{\sigma} \paren{t_{jk}-s}}
  \alpha_n(s) \, ds \\
  &= \sum_{\nu=\max(0,M-j)}^{M-2} \int_{(j-M+\nu) \Delta t}^{(j-M+\nu+1)
  \Delta t} K_{mn}\paren{\frac{c}{\sigma} \paren{t_{jk}-s}}
  \alpha_n(s) \, ds \\
  &= \sum_{\nu=\max(0,M-j)}^{M-2} \sum_{l=0}^{q-1}
  \wh{\alpha}_{n,j-M+\nu+1,l} \int_{(j-M+\nu) \Delta t}^{(j-M+\nu+1)
  \Delta t} K_{mn}\paren{\frac{c}{\sigma} \paren{t_{jk}-s}}
  P_l^{j-M+\nu+1}(s) \, ds \\
  &= \sum_{\nu=\max(0,M-j)}^{M-2} \sum_{l=0}^{q-1}
  \wh{\alpha}_{n,j-M+\nu+1,l} \int_{0}^{\Delta t}
  K_{mn}\paren{\frac{c}{\sigma}\paren{(M-\nu-1) \Delta t + \tau_k -s}}
  P_l(s) \, ds \\
  &= \sum_{\nu=\max(0,M-j)}^{M-2} \sum_{l=0}^{q-1}
  \wh{\mathcal{L}}_{m,n,k,l,\nu} \wh{\alpha}_{n,j-M+\nu+1,l} \\
  &= \sum_{\nu=\max(0,M-j)}^{M-2} \sum_{l=0}^{q-1}
  \mathcal{L}_{m,n,k,l,\nu} \alpha_{n,j-M+\nu+1,l}
\end{align*}
\endgroup
where in the last two lines, we have again defined $\wh{\mathcal{L}}$,
and then $\mathcal{L}$, by precomposition with the discrete Legendre
transform matrix. The $M p^2 q^2$ quantities
$\mathcal{L}_{m,n,k,l,\nu}$ can be precomputed. As before, using \eqref{eq:kmndef}, this only requires computing $M p q^2$
integrals. Thus the cost of computing $L_{m,n,j,k}$ for each time step is
$\OO{M p^2 q^2}$.

\subsection{The history term} \label{sec:history}

A naive treatment of the history term would simply amount to extending
the local integral back to $t = 0$ rather than $t = t_j^*$, and using
the same method. This would
require summing over the full history of the numerical solution $\alpha_{n,j,l}$
at each time step,
rather than at most the previous $M-1$ time steps, as well as
precomputing $N p q^2$ rather than $M p q^2$ integrals. We can
avoid this expense with the sum-of-exponentials history compression
technique, which has been used in a variety of contexts to compress and
efficiently update the history contribution of Volterra integral
operators \cite{alpert00,jiang04,jiang08,jiang15}. The following discussion illustrates the
technique.

We assume for now that there is a sum-of-exponentials
representation of the kernel $K_{mn}$,
\begin{equation} \label{eq:soekappa}
  K_{mn}(t) = \sum_{\mu=1}^{n_e} w_{m,n,\mu} e^{-\lambda_\mu t},
\end{equation}
valid for $\delta \leq \sqrt{2} t \leq \tmax $, with $0 < \delta < \tmax$.
Here $w_{m,n,\mu}, \lambda_\mu \in \CC$ and $\Re \lambda_\mu
\ge 0$.
We will show in Section \ref{sec:soe} that such a representation can be
constructed with $\delta = 20$, $\tmax = 10^8$, and $n_e = 67$, which is accurate to near
machine precision for all $m,n$.

Let us assume $T \leq \sigma \tmax/(\sqrt{2} c)$, and choose $M$ such that $(M-1) \Delta t \ge \sigma
\delta/(\sqrt{2} c)$. If $j
\le M$, we have $t_j^* = 0$ and $H_{m,n,j,k} = 0$. Otherwise, we have $j
\geq M+1$, so that
\begin{equation} \label{eq:hist1}
\begin{aligned}
  H_{m,n,j,k} &= \int_0^{(j-M) \Delta t} 
  K_{mn}\paren{\frac{c}{\sigma} \paren{t_{jk}-s}}
  \alpha_n(s) \, ds \\
  &= \sum_{\mu=1}^{n_e} w_{m,n,\mu} \int_0^{(j-M) \Delta t} 
  e^{-c \lambda_\mu (t_{jk}-s)/\sigma} \alpha_n(s) \, ds \\
  &= \sum_{\mu=1}^{n_e} w_{m,n,\mu} h_{n,j,k,\mu}
\end{aligned}
\end{equation}
where $h_{n,j,k,\mu} = \int_0^{(j-M) \Delta t}  
  e^{-c \lambda_\mu (t_{jk}-s)/\sigma} \alpha_n(s) \, ds$. Observe that
\begin{align*}
  h_{n,j,k,\mu} &= \int_0^{(j-M) \Delta t} e^{-c \lambda_\mu (t_{jk}-s)/\sigma}
\alpha_n(s) \, ds \\
  &= e^{-c \lambda_\mu \Delta t / \sigma} \int_0^{(j-1-M) \Delta t}
e^{-c \lambda_\mu (t_{(j-1)k}-s) / \sigma} \alpha_n(s) \, ds + \int_{(j-1-M)
  \Delta t}^{(j-M) \Delta t} e^{-c \lambda_\mu (t_{jk}-s)/\sigma}
\alpha_n(s) \, ds \\
  &= e^{-c \lambda_\mu \Delta t / \sigma} h_{n,j-1,k,\mu} +
  \int_{(j-1-M) \Delta t}^{(j-M) \Delta t} e^{-c \lambda_\mu
  (t_{jk}-s) / \sigma}
\alpha_n(s) \, ds.
\end{align*}
This is a recurrence for $h_{n,j,k,\mu}$. To update it from one time
step to the next, we multiply by a damping factor and add a local update
integral. For the local update integral, we write
\begin{align*}
  \int_{(j-1-M) \Delta t}^{(j-M) \Delta t} e^{-c \lambda_\mu
  (t_{jk}-s) / \sigma} \alpha_n(s) \, ds &= \sum_{l=0}^{q-1} \wh{\alpha}_{n,j-M,l} \int_{(j-1-M)
  \Delta t}^{(j-M) \Delta t} e^{-c \lambda_\mu (t_{jk}-s) / \sigma}
  P_l^{j-M}(s) \, ds \\
  &= \sum_{l=0}^{q-1} \wh{\alpha}_{n,j-M,l} \int_0^{\Delta t} 
  e^{-c \lambda_\mu (M \Delta t + \tau_k -s) / \sigma}
  P_l(s) \, ds \\
  &= \sum_{l=0}^{q-1} \wh{\mathcal{H}}_{k,l,\mu} \wh{\alpha}_{n,j-M,l} \\
  &= \sum_{l=0}^{q-1} \mathcal{H}_{k,l,\mu} \alpha_{n,j-M,l},
\end{align*}
where the second to last line defines $\wh{\mathcal{H}}_{k,l,\mu}$,
and $\mathcal{H}_{k,l,\mu}$ is again obtained from
$\wh{\mathcal{H}}_{k,l,\mu}$ by precomposition with the discrete
Legendre transform matrix. The $n_e q^2$ quantities
$\mathcal{H}_{k,l,\mu}$ can be precomputed. We obtain
\begin{equation} \label{eq:histupdate}
  h_{n,j,k,\mu} = e^{-c \lambda_\mu \Delta t / \sigma} h_{n,j-1,k,\mu} +
\sum_{l=0}^{q-1} \mathcal{H}_{k,l,\mu} \alpha_{n,j-M,l}
\end{equation}
which, combined with \eqref{eq:hist1}, completes our treatment of the history
term. The cost of updating $h_{n,j,k,\mu}$ at each time step using \eqref{eq:histupdate}
is $\OO{p q^2 n_e}$, and the cost of computing $H_{m,n,j,k}$ from these
values is $\OO{p^2 q n_e}$. For comparison, the cost of computing
$H_{m,n,j,k}$ directly at each time step, using the same method as we
use for the local term, would be $\OO{p^2 q^2 N}$, in addition to the
significantly larger precomputation cost.

\subsection{Summary of the time-stepping procedure and computational
complexity}

We can now summarize the full solver. We first 
precompute and store the quantities
$\mathcal{C}_{m,n,k,l}$, $\mathcal{L}_{m,n,k,l,\nu}$, and
$\mathcal{H}_{k,l,\mu}$. Now
let
\[b_{m,j,k} = - \frac{g^2}{2\pi c} \sum_{n=1}^p \paren{L_{m,n,j,k} +
H_{m,n,j,k}} + f_m(t_{jk})\]
be the right hand side of \eqref{eq:viedisc}. Sections \ref{sec:local}
and \ref{sec:history} describe how to compute $b_{m,j,k}$ at each
time step using the precomputed arrays and the values of the solution at
the previous $M$ time steps. Using this and
\eqref{eq:cmat}, we can write the discretized VIE \eqref{eq:viedisc} as
\[\alpha_{m,j,k} + \frac{g^2}{2\pi c} \sum_{n=1}^p \sum_{l=0}^q
\mathcal{C}_{m,n,k,l} \alpha_{n,j,l} = b_{m,j,k}.\]
To take the $j$th time step, we solve this $pq \times pq$ linear system.
The system matrix, with entries $\delta_{mn}\delta_{kl} + \mathcal{C}_{m,n,k,l}$, can be
formed and $LU$-factorized as a precomputation.

The cost of computing $b_{m,j,k}$
at each time step is $\OO{p^2 q^2 M + (p^2 q + p q^2) n_e}$, ignoring
the evaluation of $f_m(t)$. The
cost of solving the linear system by backward substitution is just $\OO{p^2
q^2}$.
Let us write the computational complexity in terms of the number of
time steps, $N$. There are two $N \to \infty$ regimes: $T$ fixed,
$\Delta t \to 0$, and $\Delta t$ fixed, $T \to \infty$. In practice, using
high-order time-stepping, convergence with respect to $\Delta t$ is
rapid, and the limit
$\Delta t \to 0$ is unimportant; see Figures \ref{fig:ex1_err} and
\ref{fig:ex2_err} in Section \ref{sec:results}.
With $\Delta t$ fixed and $T \to \infty$, $M$ is fixed, and the computational complexity
is $\OO{N n_e}$. $n_e$ in turn depends on $\tmax$, and in particular, as
we will discuss later, grows like $\OO{\log \tmax}$. To ensure $T = N \Delta
t \leq \tmax/\sqrt{2}$, then, we have $n_e = \OO{\log N}$, giving
overall $\OO{N \log N}$ computational complexity.

\section{Representation and evaluation of kernels} \label{sec:kernels}

We have seen that building the arrays $\mathcal{C}_{m,n,k,l}$ and
$\mathcal{L}_{m,n,k,l,\nu}$ requires computing integrals against the
kernels $k_n$. In particular, if we use standard integration
routines, we require a method of evaluating those
kernels for all $t \geq 0$. Furthermore, evolving the history term requires a
sum-of-exponentials representation \eqref{eq:soekappa} of $K_{mn}$, valid
for sufficiently large times.

We will accomplish both objectives by using
an efficient representation of the kernel $j_n$, defined by
\eqref{eq:jndef}. First, we will obtain a sum-of-exponentials
representation of $j_n$, valid when $t>\delta = 20$, and use it to
obtain similar representations for $k_n$ and hence $K_{mn}$. This also
solves the problem of evaluating $k_n$ for sufficiently large $t$. Then we will obtain
Chebyshev expansions of $j_n$ and thereby of $k_n$ valid for $t\leq\delta$.

\begin{figure}[t]
  \centering
  \begin{subfigure}[t]{.3\textwidth}
    \centering
    \includegraphics[width=\linewidth]{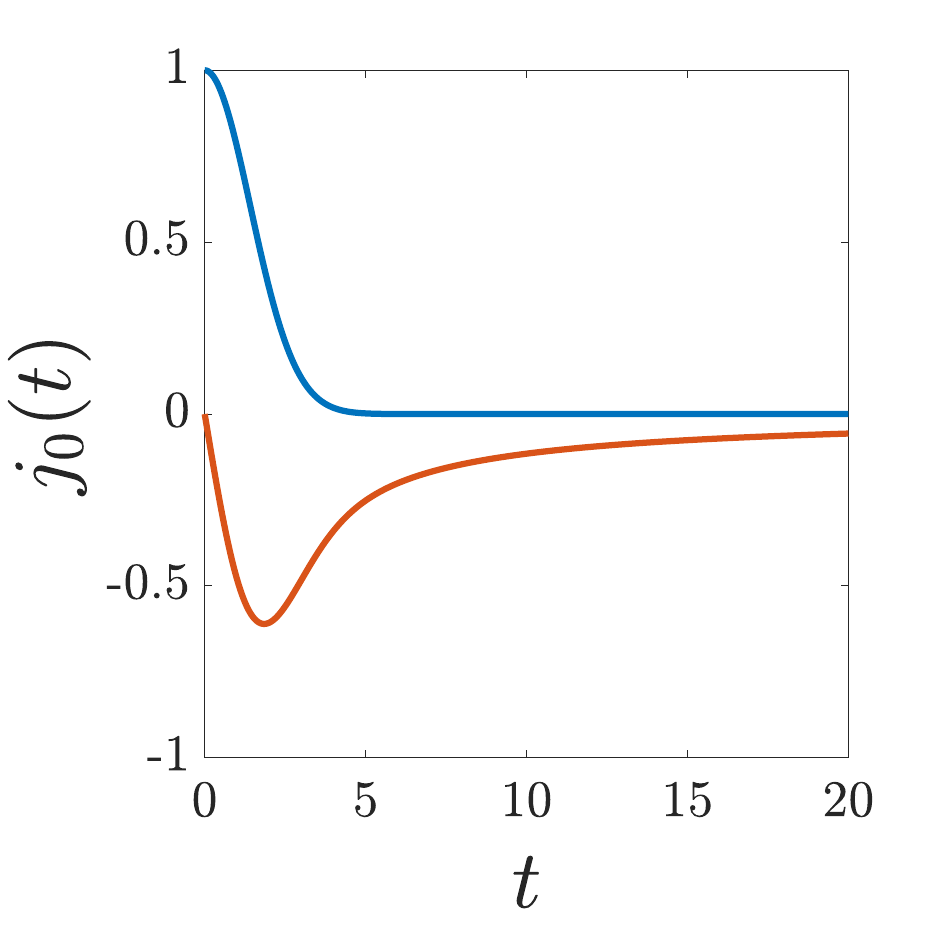}
    \label{fig:j0_short}
  \end{subfigure}
  \begin{subfigure}[t]{.3\textwidth}
    \centering
    \includegraphics[width=\linewidth]{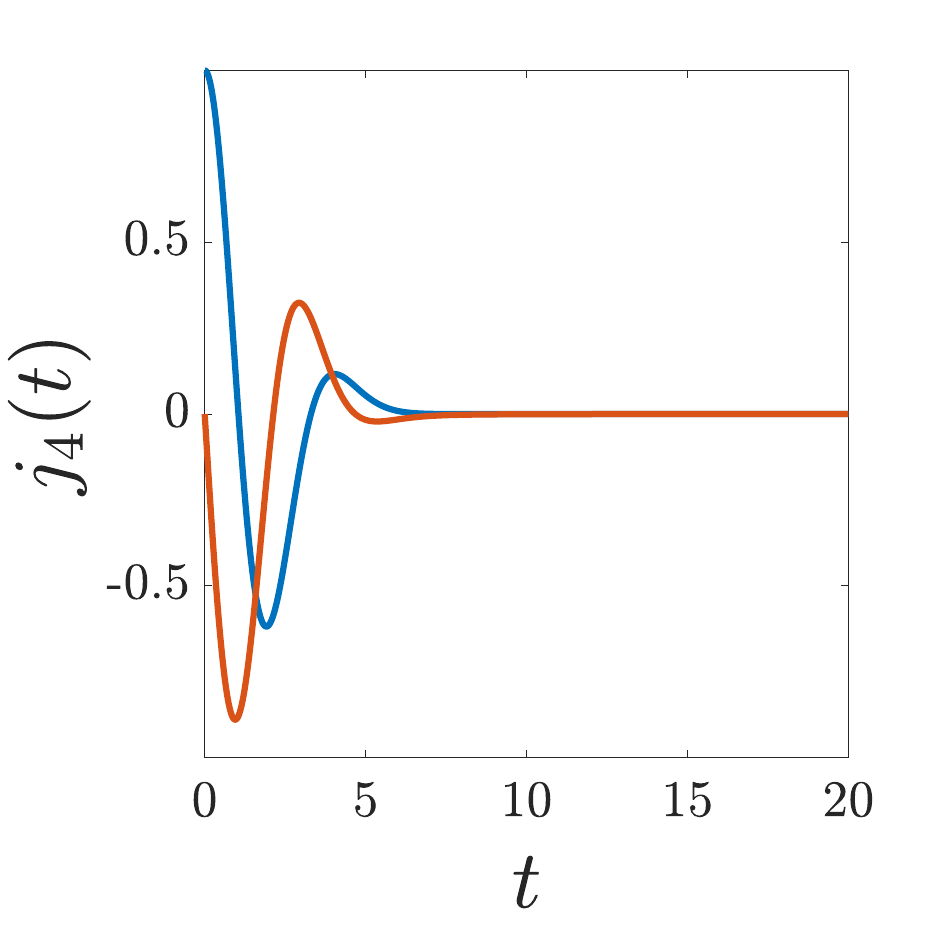}
    \label{fig:j4_short}
  \end{subfigure}
  \begin{subfigure}[t]{.3\textwidth}
    \centering
    \includegraphics[width=\linewidth]{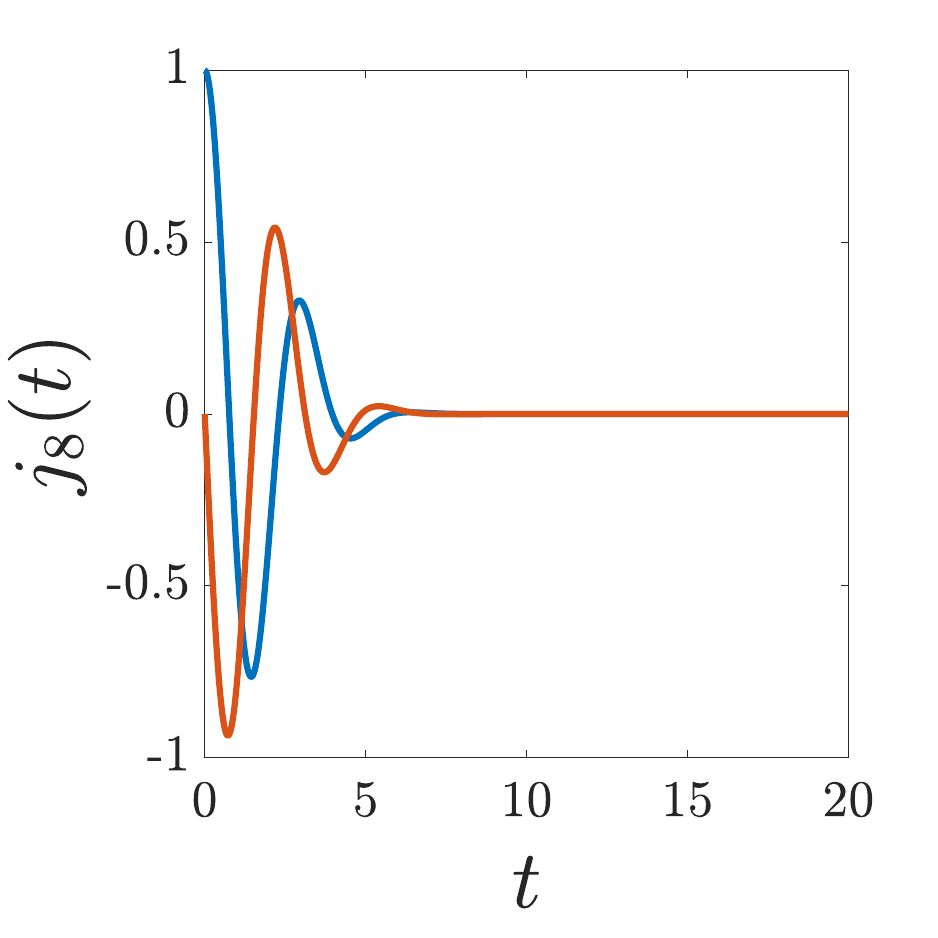}
    \label{fig:j8_short}
  \end{subfigure}

  \begin{subfigure}[t]{.3\textwidth}
    \centering
    \includegraphics[width=\linewidth]{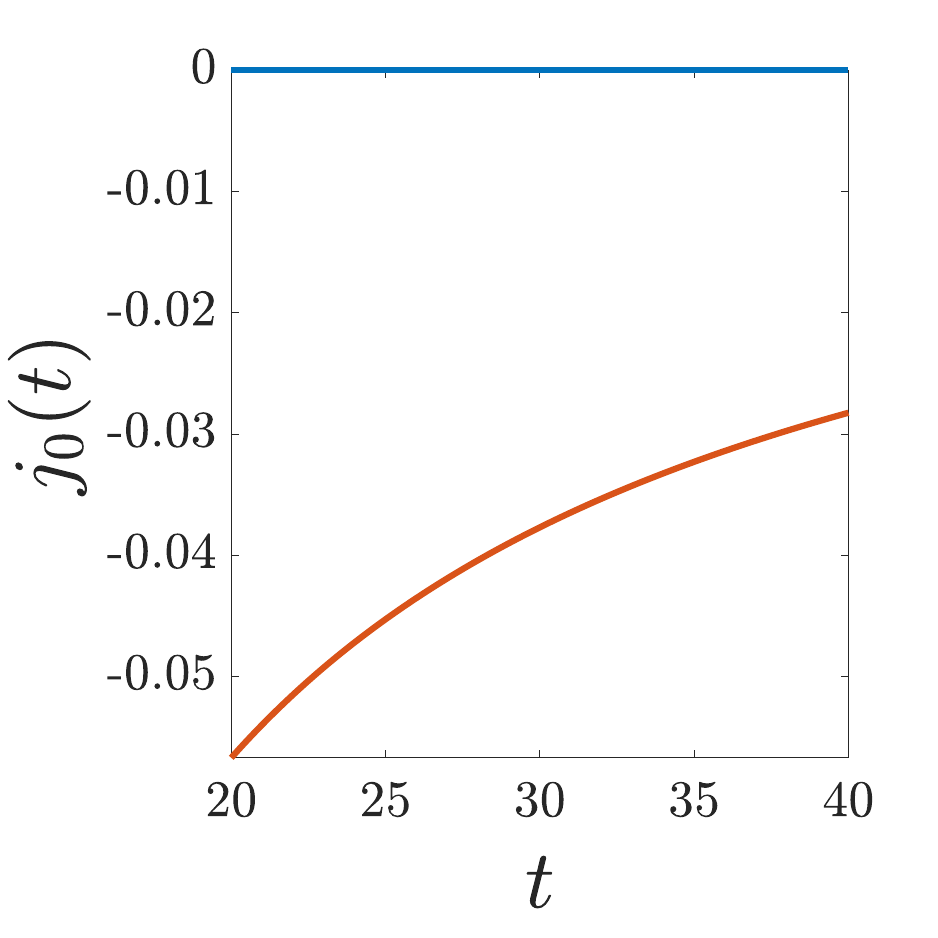}
    \label{fig:j0_long}
  \end{subfigure}
  \begin{subfigure}[t]{.3\textwidth}
    \centering
    \includegraphics[width=\linewidth]{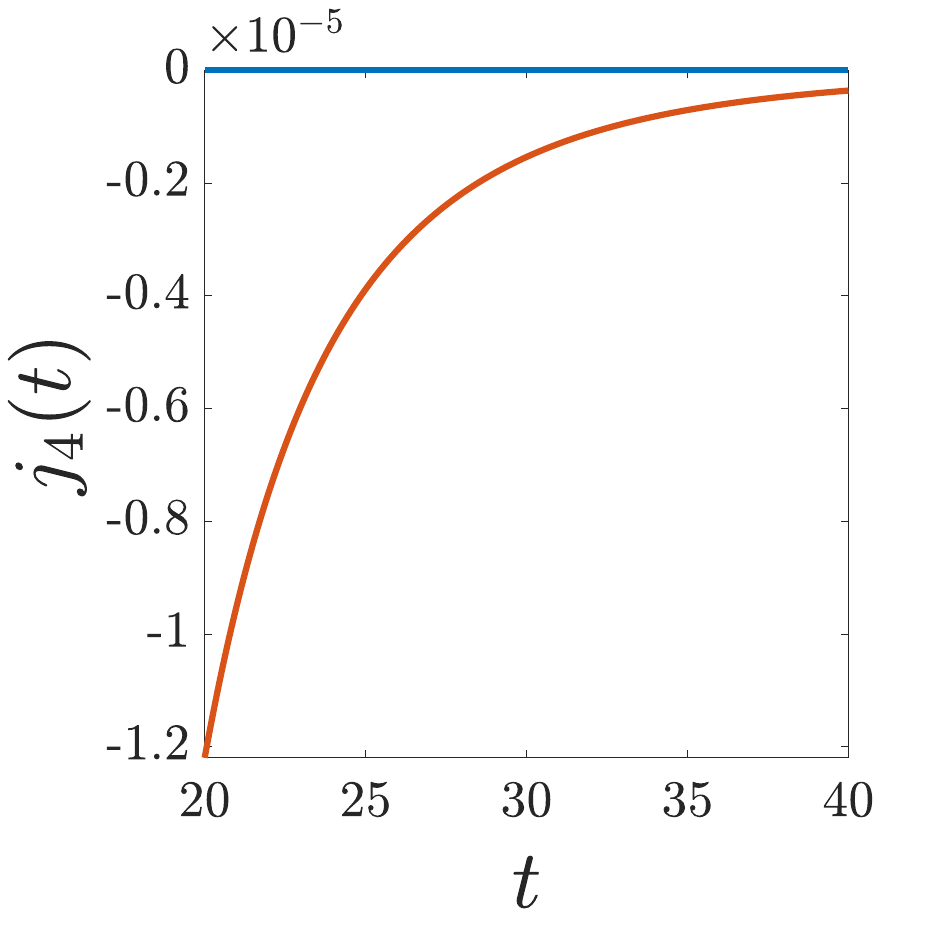}
    \label{fig:j4_long}
  \end{subfigure}
  \begin{subfigure}[t]{.3\textwidth}
    \centering
    \includegraphics[width=\linewidth]{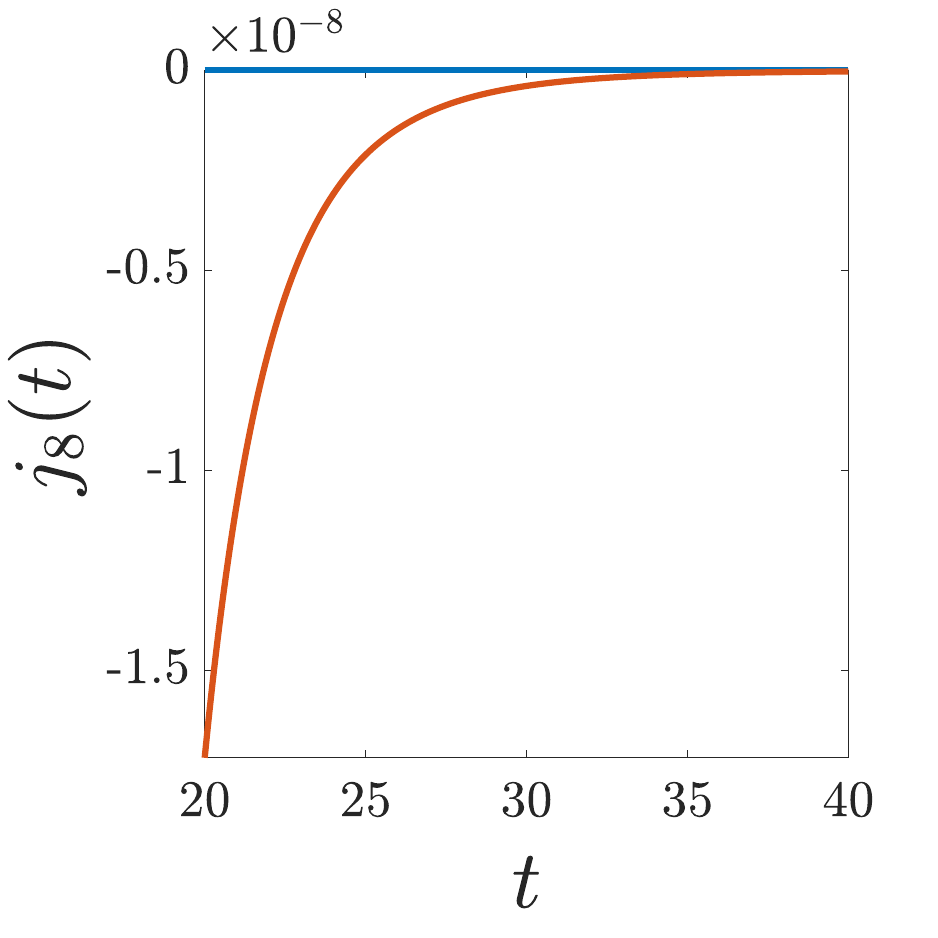}
    \label{fig:j8_long}
  \end{subfigure}

  \caption{The first row shows $j_n(t)$ for $n = 0,4,8$, and $t$
  in the small-time interval $[0,20]$.
  $\Re{j_n}(t)$ is
  indicated by the blue curve, and $\Im{j_n}(t)$ by the red curve. In this
  interval, we represent each $j_n$ by a Chebyshev expansion. The second
  row shows the same kernels for $t > 20$, where we represent
  them by sum-of-exponentials expansions.}
\label{fig:kmn}
\end{figure}

Figure \ref{fig:kmn} shows representative examples of the kernels
$j_n(t)$ for $t \in [0,20]$ and $t > 20$.

\subsection{Sum-of-exponentials representation for large times}
\label{sec:soe}

We start by constructing a sum-of-exponentials representation of
$j_n(t)$ valid for $t > \delta$, for some $\delta > 0$ to be determined. We note that a
sum-of-exponentials representation
\eqref{eq:soekappa} of $K_{mn}$ can then be obtained from
\eqref{eq:kndef} and \eqref{eq:kmndef}. In particular,
suppose the representation
\begin{equation} \label{eq:ksoe}
  j_n(t) = \sum_{\mu=1}^{n_e} \wt{w}_{n,\mu} e^{-\wt{\lambda}_\mu t}
\end{equation}
is valid for $t > \delta$, for $\wt{w}_{n,\mu} \in
\CC$ and $\wt{\lambda}_\mu > 0$. Then for $t
\geq \delta/\sqrt{2}$, when $m+n$ is even, we have
\begin{align*}
  K_{mn}(t) &= K_{mn}\paren{\delta/\sqrt{2}} + (-1)^m (-i)^{m+n}
  \frac{\Gamma\paren{\frac{m+n+1}{2}}}{\sqrt{\frac{m! n!}{2}}}
  \int_{\delta/\sqrt{2}}^t e^{i \frac{\Omega \sigma}{c} s}
  j_{m+n}(\sqrt{2} s) \, ds \\
  &= K_{mn}\paren{\delta/\sqrt{2}} + (-1)^m (-i)^{m+n}
  \frac{\Gamma\paren{\frac{m+n+1}{2}}}{\sqrt{\frac{m! n!}{2}}} \sum_{\mu=1}^{n_e} \wt{w}_{n,\mu} 
  \int_{\delta/\sqrt{2}}^t e^{\paren{i \frac{\Omega \sigma}{c} - \sqrt{2}
  \wt{\lambda}_\mu} s} \, ds \\
  &= \begin{multlined}[t][.74\textwidth]
    K_{mn}\paren{\delta/\sqrt{2}} + (-1)^m (-i)^{m+n}
  \frac{\Gamma\paren{\frac{m+n+1}{2}}}{\sqrt{\frac{m! n!}{2}}} \sum_{\mu=1}^{n_e} \frac{\wt{w}_{n,\mu}}{i
  \Omega\sigma/c - \sqrt{2} \wt{\lambda}_\mu} \\
    \times \paren{e^{\paren{i
  \frac{\Omega \sigma}{c} - \sqrt{2} \wt{\lambda}_\mu} t} - e^{\paren{i
  \frac{\Omega \sigma}{c} - \sqrt{2} \wt{\lambda}_\mu} \delta/\sqrt{2}}}
  \end{multlined}
  \\
  &= \sum_{\mu=1}^{n_e+1} w_{m,n,\mu} e^{-\lambda_\mu t}
\end{align*}
with
\[w_{m,n,\mu} = 
\begin{cases}
  (-1)^m (-i)^{m+n} \frac{\Gamma\paren{\frac{m+n+1}{2}}}{\sqrt{\frac{m!
  n!}{2}}}
  \frac{\wt{w}_{n,\mu}}{i \Omega \sigma / c - \sqrt{2} \wt{\lambda}_\mu}
  & \text{if } 1 \leq \mu \leq
  n_e \\
  K_{mn}\paren{\delta/\sqrt{2}} - (-1)^m (-i)^{m+n} \frac{\Gamma\paren{\frac{m+n+1}{2}}}{\sqrt{\frac{m!
  n!}{2}}} \sum\limits_{\nu=1}^{n_e} \frac{\wt{w}_{n,\nu}}{i
  \Omega \sigma/c - \sqrt{2} \wt{\lambda}_\nu} e^{\paren{i
  \frac{\Omega \sigma}{c} - \sqrt{2} \wt{\lambda}_\nu} \delta/\sqrt{2}}
  & \text{if } \mu = n_e+1
\end{cases}
\]
and 
\[\lambda_\mu = 
\begin{cases}
  \sqrt{2} \wt{\lambda}_\mu - i \Omega \sigma / c & \text{if } 1 \leq \mu \leq n_e \\
  0 & \text{if } \mu = n_e+1.
\end{cases}
\]
Redefining $n_e \gets n_e+1$ gives the desired representation \eqref{eq:soekappa}. We focus then on
the construction of \eqref{eq:ksoe}.

\begin{figure}[t]
  \centering
  \includegraphics[width=0.5\linewidth]{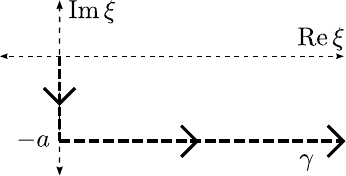}
  \caption{Contour of integration used for $J_{mn}(t)$ to obtain
  sum-of-exponentials representation. For sufficiently large $a$ and
  $\delta$, the contribution from the
  horizontal part of the contour is negligible when $t > \delta$. The
  sum-of-exponentials representation is obtained by applying a quadrature rule
  to the vertical part.}
  \label{fig:contour}
\end{figure}

We begin by deforming the integral defining $j_n$ in \eqref{eq:jndef}
from the interval $[0,\infty)$ to the contour $\gamma$ shown in Figure \ref{fig:contour}. That is, we have
\begin{align*}
  j_n(t) &= \frac{2}{\Gamma\paren{\frac{n+1}{2}}} \paren{(-i)^{n+1} \int_0^a
  \eta^n e^{\eta^2 - \eta t} \, d\eta + \int_0^\infty
  (\eta-ia)^n e^{-(\eta-ia)^2 - i (\eta-i a) t} \, d\eta} \\
  &\equiv j_n^{(1)}(t) + j_n^{(2)}(t).
\end{align*}
We show in Appendix \ref{sec:jn2} that
\[\abs{j_n^{(2)}(t)} \leq 14 e^{2a^2 - at}.\]
If we take $a = \delta/4$, then $\abs{j_n^{(2)}(t)} \leq 14
e^{-\delta^2/8}$ when $t \geq \delta$. Thus, to ensure $\abs{j_n^{(2)}(t)} < \varepsilon$ for all $t \geq
\delta$, we can take $\delta > \sqrt{8 \log(14/\varepsilon)}$. If $\varepsilon$ is the double machine
precision, then taking $\delta = 20$, $a = 5$ is sufficient to neglect
$\abs{j_n^{(2)}(t)}$.

As a consequence, if we can find a quadrature rule $\{\omega_{n,\mu},\wt{\lambda}_\mu\}_{\mu=1}^{n_e}$ so
that
\[j_n^{(1)}(t) = \frac{2  (-i)^{n+1}}{\Gamma\paren{\frac{n+1}{2}}}
\int_0^a \eta^{n} e^{\eta^2-\eta t} \, d\eta \approx
\sum_{\mu=1}^{n_e} \omega_{n,\mu} \wt{\lambda}_\mu^n e^{\wt{\lambda}_\mu^2-\wt{\lambda}_\mu t}\]
holds to high accuracy for all $t \geq \delta$,
then this gives \eqref{eq:ksoe} with
\[\wt{w}_{n,\mu} = \omega_{n,\mu} \wt{\lambda}_\mu^n e^{\wt{\lambda}_\mu^2}.\]
When $n > 23$, $\abs{j_n(t)}$ is below the double
machine precision for all $t > 20$. We therefore only need quadratures for the above integrals which are valid for $n
= 0,\ldots,23$; we can simply take $\omega_{n,\mu} = 0$ for $n \geq 24$.

The method of generalized Gaussian quadrature can be used to find such a
quadrature rule \cite{ma96}. Given a family of functions -- in this case, the
functions $\eta^n e^{\eta^2-\eta t}$, $\eta \in [0,5]$, for $n = 0,\ldots,23$
and $\delta < t < \tmax$ -- this method uses a nonlinear optimization process to
determine a minimal set of quadrature nodes and weights sufficient to
integrate all functions in the family to near machine precision. An
upper bound on the number of quadrature nodes required can be given in
terms of the numerical rank of the family of functions. It is
straightforward to adapt the proof given in Ref. \cite[Lemma 4.4]{gimbutas20} for the
case of a family of decaying exponentials to the present setting.
Briefly, the
proof works by 1) rescaling the interval to $[0,1]$; 2) discretizing
$[0,1]$ by a composite Chebyshev grid with nodes exponentially
clustered at the origin; and 3) using standard error estimates for Chebyshev
interpolation to show that the resulting piecewise polynomial
approximation is uniformly accurate for all functions in the family.
This argument shows that the numerical rank of the family scales as $\OO{\log
\paren{\tmax/\delta}}$.
In practice, we simply take $\tmax =
10\,000\,000$, several orders of magnitude larger than is needed for the
examples shown in this article, and obtain a quadrature rule of $n_e =
67$ nodes and weights.

\subsection{Chebyshev representation for small times}

We next consider the evaluation of $k_n(t)$ for $t \leq
\delta/\sqrt{2}$. First, we can evaluate each $j_n(t)$ at Chebyshev nodes on
$[0,\delta]$ using adaptive integration.
$j_n(t)$ is an entire function, so its Chebyshev
interpolant converges super-exponentially \cite{trefethen19}. A moderate number of
Chebyshev nodes are therefore sufficient to represent the function on the
full interval to near machine precision by its interpolant at these nodes; see Figure \ref{fig:kmn} for plots of some $j_n(t)$ on
$[0,\delta]$. The samples at Chebyshev nodes
can be computed once and stored. $j_n(t)$ can then be evaluated at any
$t \in [0,\delta]$ by barycentric interpolation
\cite{trefethen19,berrut04,higham04}.

Given $\Omega$, $c$, and $\sigma$, samples of the
integrand of $k_n(t)$ in \eqref{eq:kndef} at Chebyshev nodes on $[0,\delta/\sqrt{2}]$ can
then be obtained by pointwise multiplication. If $\Omega \sigma/c$ is
large, then to resolve the complex exponential, $j_n$ can be evaluated
on a denser Chebyshev grid. Accurate samples of $k_n(t)$ at the same Chebyshev nodes can then be
obtained by spectral integration \cite{greengard91}, and as
before, can be used to represent $k_n(t)$ on $[0,\delta/\sqrt{2}]$ by
barycentric interpolation.

\section{Numerical results} \label{sec:results}

We demonstrate the solver using two examples. In the first, we
place an atom in its excited state and observe its decay. In the second, 
we excite the atom with a wavepacket.

\subsection{Example 1: decay of an excited atom}

The first example is characterized by the initial condition
\[a(x,0) = 1, \quad u(x,0) = 0,\]
which corresponds to taking $\alpha_m(0) = \delta_{0n}$ and $U_m(t)
= 0$ in \eqref{eq:vie}. We take $c = \Omega = 1$ and $\sigma =
0.1$.

\begin{figure}[t]
  \centering
  \includegraphics[width=\linewidth]{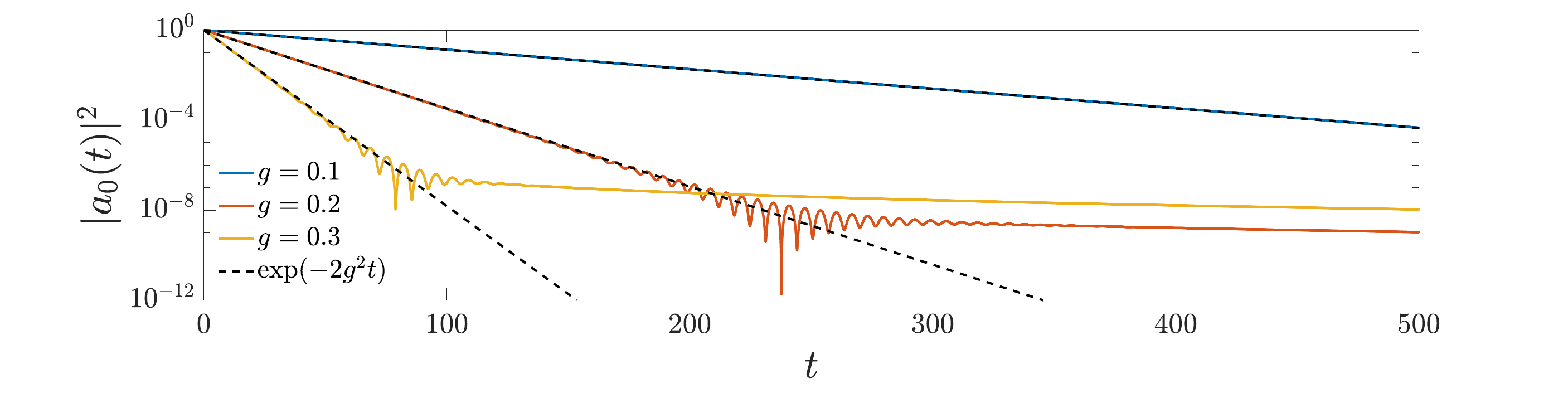}
  \caption{$\abs{a_0(t)}^2$ for the first example with $p = 1$ and
  several choices of $g$, along with the expected initial decay curve.}
  \label{fig:a_onemode}
\end{figure}

We first represent the solution using only a single Hermite polynomial, $p =
1$. Figure  \ref{fig:a_onemode} shows $\abs{a_0(t)}^2$ for $g = 0.1, 0.2,
0.3$. The solutions are characterized by an initial exponential decay
regime, with the decay rate determined by $g$, followed by a tail of
significantly slower decay. The plots indicate close agreement with the 
standard Wigner-Weisskopf estimate $\abs{a_0(t)}^2 \approx e^{-2 g^2 t}$
of the initial decay rate; for a derivation in the three-dimensional
case, which is straightforwardly adapted to the one-dimensional case, we refer to \cite{kraisler22}.

\begin{figure}[t]
  \centering
  \includegraphics[width=\linewidth]{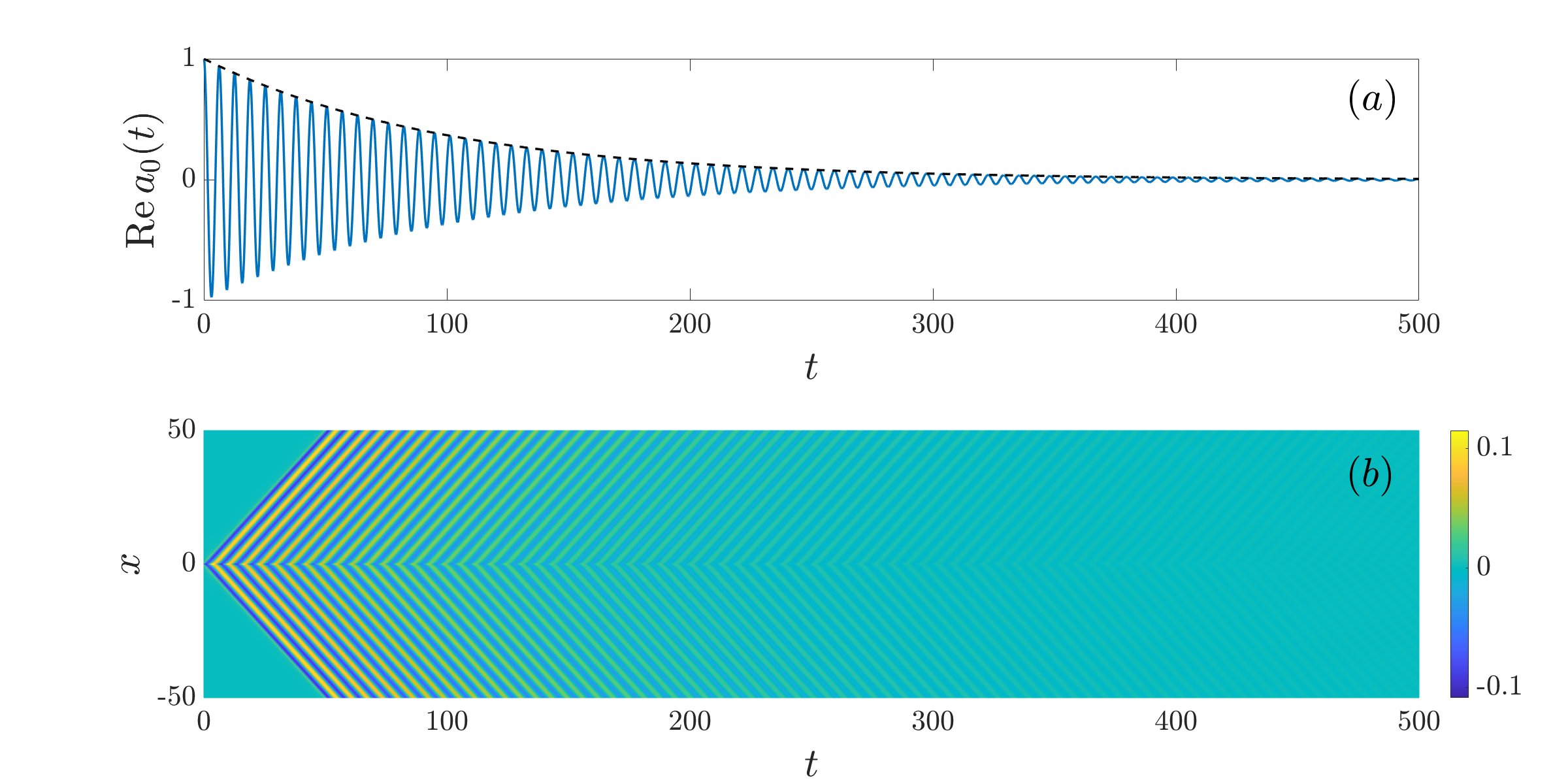}
  \caption{(a) $\Re a_0(t)$ and (b) $\Re u(x,t)$ for the first example
  with $p = 1$ and $g = 0.1$. In (a), the black dashed line is the
  curve $\exp(-g^2 t)$.}
  \label{fig:u_onemode}
\end{figure}

Figure \ref{fig:u_onemode} shows $\Re a_0(t)$ and $\Re u(x,t)$ for $g =
0.1$. As the atom amplitude decays, it acts as a source for the photon
field, which resembles a wave of speed $c$ radiating from the origin.
We note that the photon amplitude is not
identically zero outside of the light cone associated with
speed $c$. Rather, as a result of nonlocal effects arising from the
fractional Laplacian term of \eqref{eq:system}, it decays
algebraically outside of the light cone.

\begin{figure}[t]
  \centering
  \includegraphics[width=\linewidth]{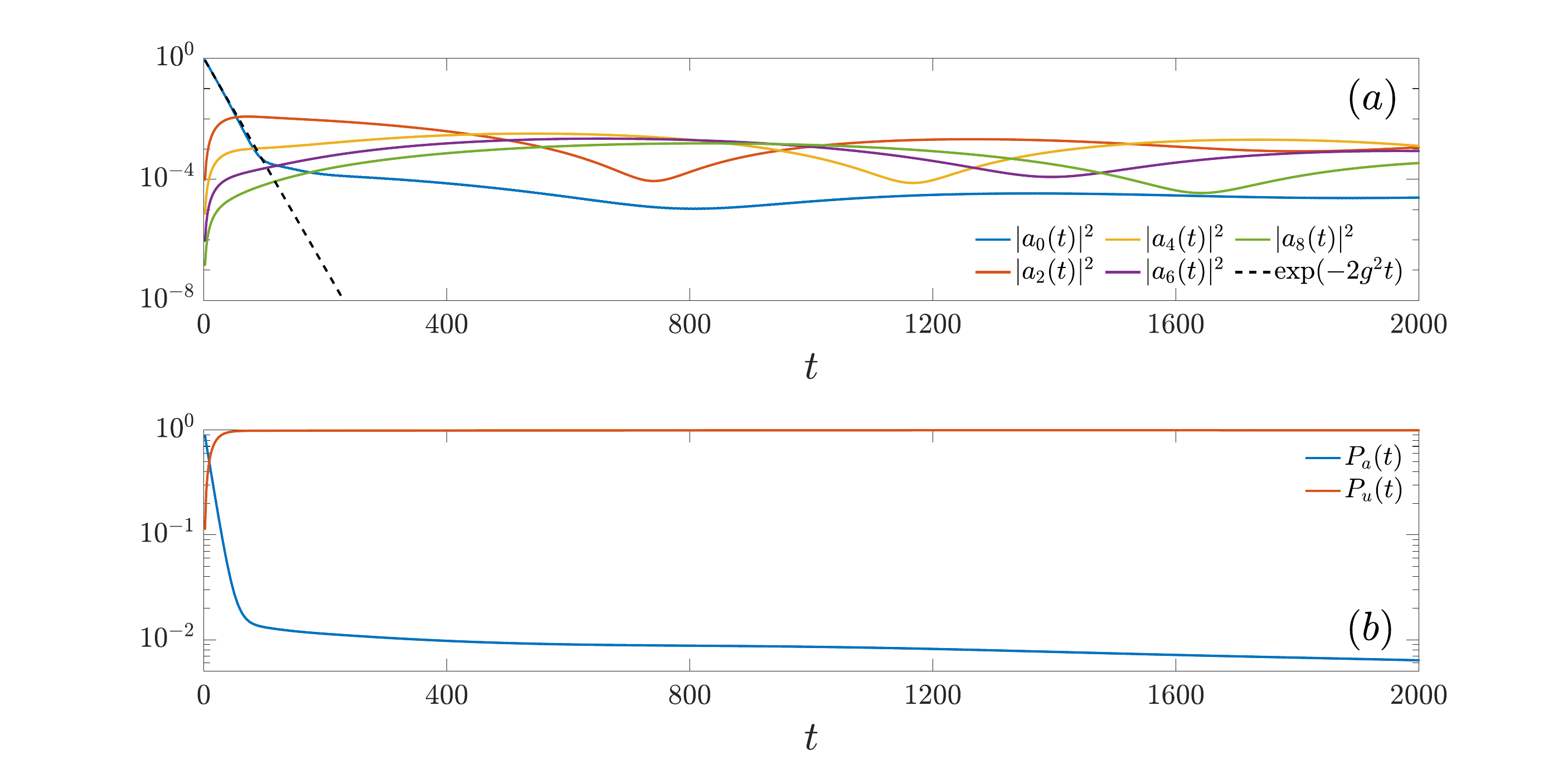}
  \caption{(a) $\abs{a_n(t)}^2$ with $n = 0,2,\ldots,8$ for the first
  example with $p = 40$ and $g = 0.2$. (b) The total probability
  associated with the atom and the photon field.}
  \label{fig:a_infmode}
\end{figure}

We next consider the limit of a large number of Hermite polynomials, $p \to
\infty$. To do so, we increase $p$ until the first five non-zero
coefficients $a_n$ are converged to high accuracy (we note that, as a
result of the symmetry of the system about the origin, the odd
coefficients are identically zero). $p = 40$ was
sufficient for the simulations considered here. We fix $g =
0.2$, so our results can be compared with the red curve in Figure
\ref{fig:a_onemode}. Figure
\ref{fig:a_infmode}a shows $\abs{a_n(t)}^2$ for the first five even
Hermite polynomials. At very short times, the $n = 0$
coefficient decays with the same rate as in the $p = 1$ case. However, the
rapid decay regime ends sooner, and gives way to
complicated, long-lived dynamics among the coefficients of the various
Hermite polynomials.

Another perspective is given by Figure \ref{fig:a_infmode}b, which shows
the total probability associated with the atom and photon
amplitudes. It can also be compared directly with the red curve in
Figure \ref{fig:a_onemode}, since in the $p = 1$ case $\abs{a_0(t)}^2 =
P_a(t)$. Compared with the $p = 1$ case, in the $p = 40$ case the
atom dissipates much less of its probability mass into the photon field. 

Evidently, allowing multiple Hermite coefficients gives rise to a trapping
effect, whereby some portion of the probability associated with the zero
coefficient
remains trapped in higher-order modes rather than being radiated into
the continuum through the photon field. As the solution
evolves in time, higher and higher-order coefficients become activated,
and the total probability associated with the atomic amplitude decays
exceptionally slowly, if at all.

The plot of the photon amplitude for the $p = 40$
case is qualitatively similar
to that appearing in Figure \ref{fig:u_onemode} for the $p = 1$ case, so
we do not plot it. The main difference is that in the $p = 40$ case,
the atom remains a longer-lived source of larger magnitude for the photon field.

\begin{figure}[t]
  \centering
  \includegraphics[width=\linewidth]{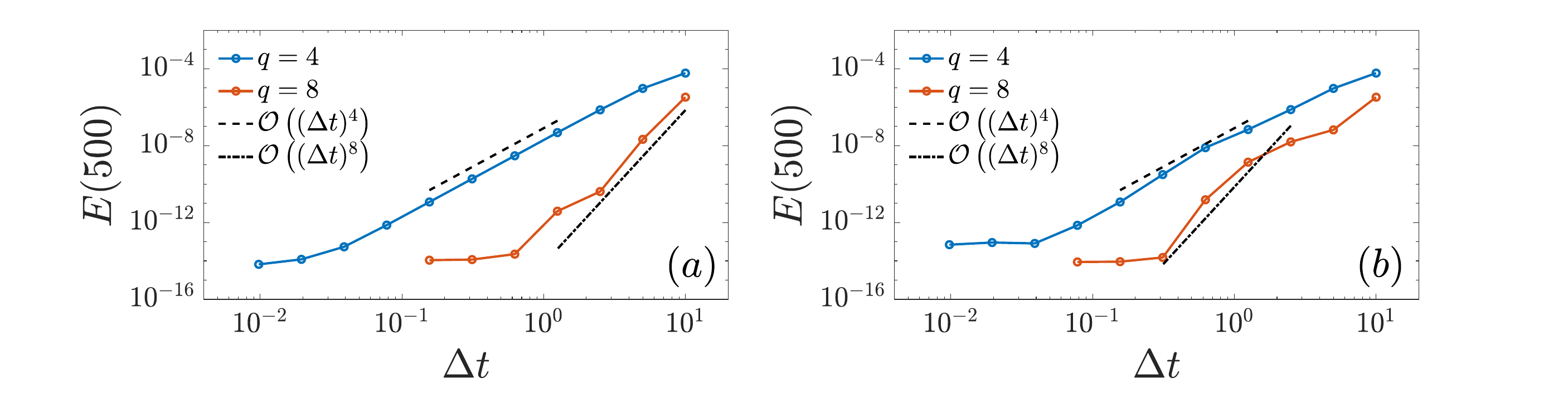}
  \caption{Error $E(t)$ for the first example with $g = 0.2$ for (a) $p
  = 1$ and (b) $p = 40$, using fourth and eighth-order time-stepping.}
  \label{fig:ex1_err}
\end{figure}

We next verify the order of accuracy of the time-stepping
algorithm by measuring the error
\begin{equation} \label{eq:err}
E(t) = \sqrt{\frac{1}{\sigma} \int_{-\infty}^\infty
\abs{a(x,t)-a^{\textrm{ref}}(x,t)} \rho(x/\sigma) \, dx} = \sqrt{\sum_{n=0}^{p-1}
\abs{\alpha_n(t) - \alpha_n^{\textrm{ref}}(t)}^2}
\end{equation}
against a well-converged reference solution $a^{\textrm{ref}}$ with
Hermite coefficients $\alpha_n^{\textrm{ref}}$. We take the parameters
as above with $g = 0.2$, and measure the error $E(t)$ at $t = 500$
using the fourth and eighth-order methods; $q = 4$ and $q = 8$,
respectively. Figure \ref{fig:ex1_err}a gives results for the $p = 1$ case, and
Figure \ref{fig:ex1_err}b for $p = 40$.

\subsection{Example 2: response to a photon pulse}

Our second example models the response of the atom to a
photon pulse. We take
\begin{equation} \label{eq:u0gaussian}
  u_0(x) = \paren{\frac{2}{\pi \beta^2}}^{1/4} e^{-(x-x_0)^2/\beta^2} e^{i \xi_0 x}
\end{equation}
with $x_0$, $\beta$, and $\xi_0$ the inital wavepacket center, width, and
wavenumber, respectively. The normalization ensures
$\int_{-\infty}^\infty \abs{u_0(x)}^2 \, dx = 1$. The free evolution is given by
\begingroup
\allowdisplaybreaks
\begin{align*}
  U(x,t) &= \int_{-\infty}^\infty G(x-y,t) u_0(y) \, dy \\
  &= \frac{1}{2\pi} \int_{-\infty}^\infty \wh{G}(\xi,t) e^{i \xi x}
\wh{u_0}(\xi) \, d\xi \\
  &= \frac{\sqrt{\beta}}{(2\pi)^{3/4}} \int_{-\infty}^\infty e^{-i \paren{c \abs{\xi} t - \xi x
  + x_0 (\xi-\xi_0)}}
  e^{-(\xi-\xi_0)^2 \beta^2 / 4}\, d\xi.
\end{align*}
\endgroup
If we take $\xi_0,\beta$ sufficiently large so that
$\varepsilon > \frac{\sqrt{\beta}}{(2\pi)^{3/4}} 
\int_{-\infty}^0 e^{-(\xi-\xi_0)^2 \beta^2 / 4} =
\frac{\erfc(\beta \xi_0/2)}{(8\pi \beta^2)^{1/4}}$ for some
$\varepsilon$, then
up to an error $\varepsilon$, we simply recover a translation of the initial wavepacket:
\[U(x,t) \approx \frac{\sqrt{\beta}}{(2\pi)^{3/4}} e^{i x_0 \xi_0} \int_{-\infty}^\infty
e^{i \xi \paren{x-x_0-c t}}
  e^{-(\xi-\xi_0)^2 \beta^2 / 4}\, d\xi = u_0(x-ct).\]
Thus the free evolution of a wavepacket with a sufficiently high frequency
modulation relative to its width is approximately given by 
translation at velocity $c$. We will choose $\xi_0,\beta$ so that the approximate
equality holds to machine precision -- $\xi_0 \beta \geq 12$
with $\beta \geq 1$ is
sufficient -- and for simplicity take
it to be an equality going forward.

To compute the source term in the VIE \eqref{eq:vie}, we write
\begin{equation} \label{eq:srcterm}
  \int_0^t e^{i \Omega s} \wh{u_m}(s) \, ds = \frac{1}{\sigma}
  \int_{-\infty}^\infty \rho(x/\sigma) f_m(x/\sigma) \int_0^t e^{i
  \Omega s} U(x,s) \, ds \, dx.
\end{equation}
We have
\begin{align*}
    \int_0^t e^{i \Omega s} U(x,s) \, ds &= \int_0^t e^{i \Omega s}
    u_0(x-cs) \, ds \\
    &= \paren{\frac{2}{\pi \beta^2}}^{1/4} e^{i \xi_0 x} \int_0^t e^{i (\Omega-\xi_0 c) s}
    e^{-(x-x_0-cs)^2/\beta^2} \, ds \\
    &= \begin{multlined}[t][0.765\textwidth]
      \frac{i \pi^{1/4} \sqrt{\beta}}{2^{3/4} c} e^{i \paren{\xi_0 x_0 +
  \Omega(x-x_0)/c}} e^{-\frac{\beta^2 (\Omega - \xi_0 c)^2}{4 c^2}}
    \paren{\erfi\left(\frac{\beta (\Omega - \xi_0 c)}{2c} - i
  \frac{x-x_0}{\beta}}\right. \\
    \left. - \erfi\paren{\frac{\beta (\Omega - \xi_0 c)}{2c} - i
    \frac{x-x_0-ct}{\beta}}\right).
    \end{multlined}
\end{align*}
The outer integral in \eqref{eq:srcterm} can be computed at each time step by
adaptive integration using the explicit expression for the inner
integral. To improve the efficiency, the integrals for each $m$ can
be computed simultaneously using the recurrence for the normalized Hermite functions:
\[\sqrt{\frac{m+1}{2}} f_{m+1}(x) = x f_m(x) -
\sqrt{\frac{m}{2}} f_{m-1}(x).\]

In this setup, we take $a(x,0) = 0$, with $c = \Omega = 1$ and $\sigma =
0.1$ as before, and we fix $g = 0.2$. We take $\beta = 12$ and $\mu =
-80$ in \eqref{eq:u0gaussian}, so that, to machine precision, the
wavepacket does not initially overlap with the atomic density.

\begin{figure}[t]
  \centering
  \includegraphics[width=\linewidth]{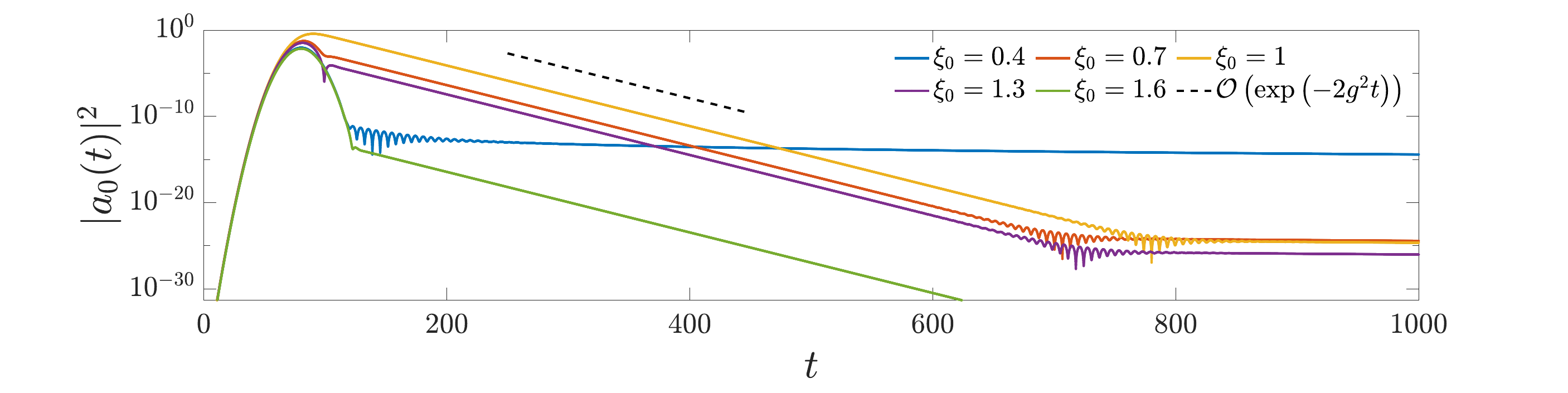}
  \caption{$\abs{a_0(t)}^2$ for the second example with $p = 1$ and
  several choices of $\xi_0$, along with an indication of the decay rate
  expected for an initially excited atom. The atomic resonance frequency
  is $\Omega = 1$.}
  \label{fig:ex2_p1_decay}
\end{figure}

We first consider the single coefficient case $p=1$. In Figure
\ref{fig:ex2_p1_decay}, we plot $\abs{a(t)}^2$ for different choices of
the wavenumber, $\xi_0 = 0.4, 0.7, 1, 1.3, 1.6$. The incoming wavepacket interacts with the atom,
increasing the magnitude of the atom amplitude, which then decays at the expected rate.
We note that in this case, the rapid decay regime continues for longer
than in the first example; a comparison can be made with the red curve in Figure
\ref{fig:a_onemode}. We also see that a wavepacket with $\xi_0 = \Omega$ --
exactly resonant with the atom -- yields the largest and most long-lived
atomic excitation. By contrast, when the modulation is chosen
off-resonance, the atomic amplitude first follows the profile of the
wavepacket-induced forcing before eventually settling into the usual
decay regime.

\begin{figure}[t]
  \centering
  \includegraphics[width=\linewidth]{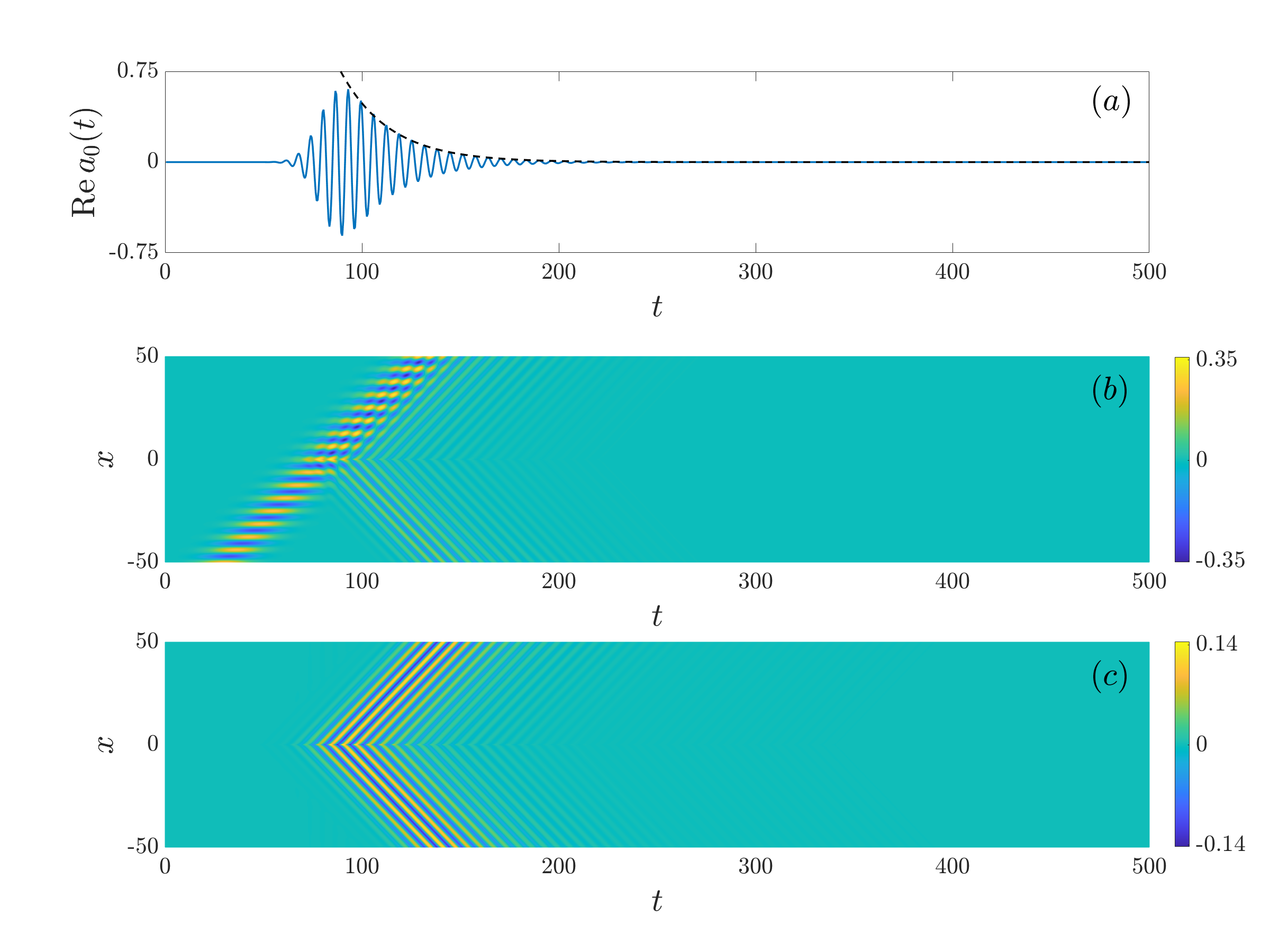}
  \caption{(a) $\Re a_0(t)$, (b) $\Re u(x,t)$, and (c)
  $\Re\paren{u(x,t)-U(x,t)}$ for the second example
  with $p = 1$ and $\xi_0 = 1$. In (a), the black dashed line is the
  curve $\exp\paren{-g^2 (t-t_0)}$ with $t_0$ chosen to follow the
  envelope of the oscillation.}
  \label{fig:ex2_p1_aandu}
\end{figure}

For the on-resonance case, $\xi_0 = 1$, Figures \ref{fig:ex2_p1_aandu}a
and \ref{fig:ex2_p1_aandu}b give plots of $\Re a(t)$ and
$\Re u(x,t)$, respectively. We see the wavepacket approach the atom
center and excite the atom, which then decays and induces its own
response in the photon amplitude, given by $u-U$. Figure
\ref{fig:ex2_p1_aandu}c gives a plot of $\Re\paren{u(x,t)-U(x,t)}$.

\begin{figure}[t]
  \centering
  \includegraphics[width=\linewidth]{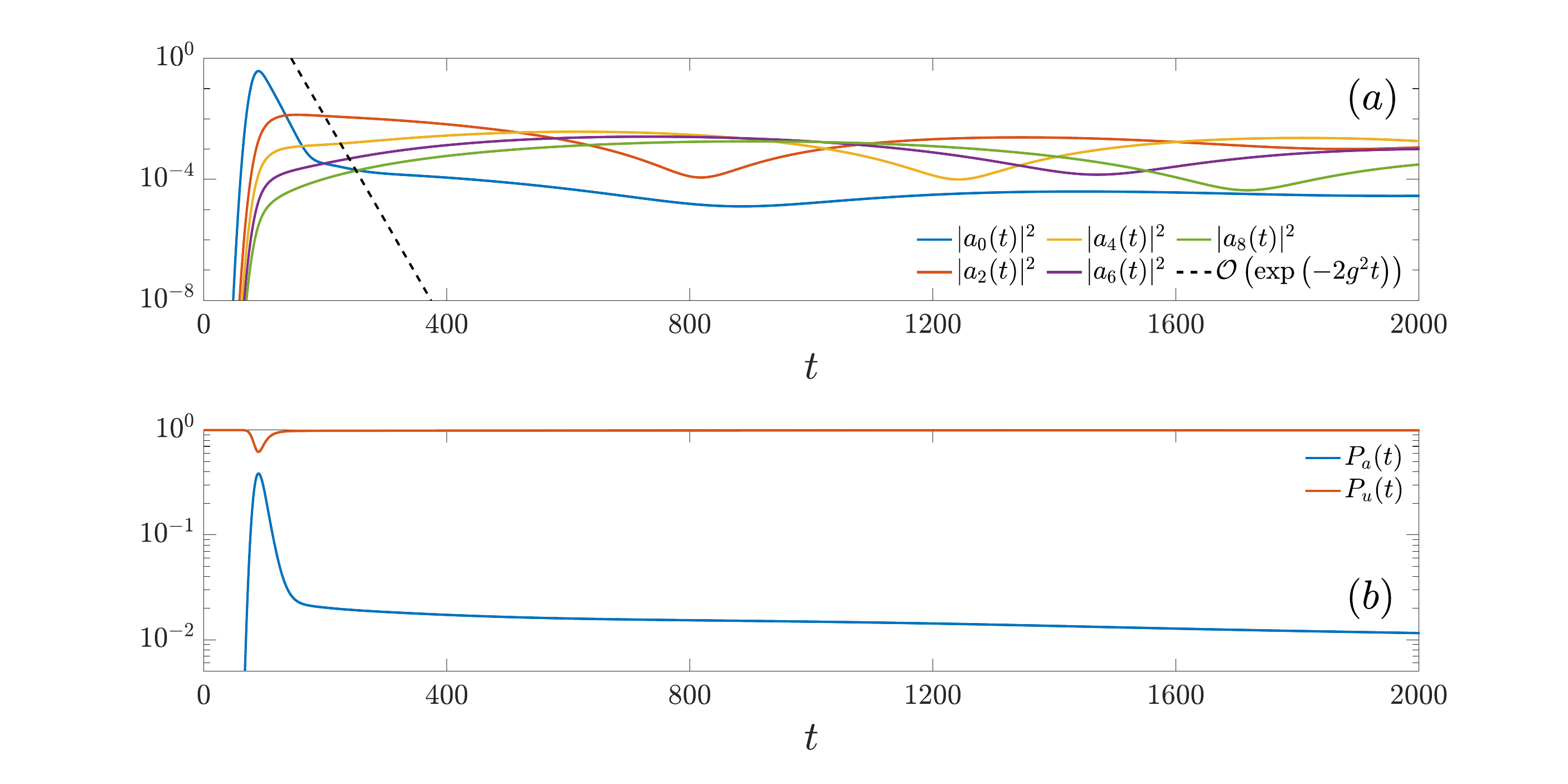}
  \caption{(a) $\abs{a_n(t)}^2$ with $n = 0,2,\ldots,8$ for the second
  example with $p = 40$ and $\xi_0 = 1$, along with an indication of
  the decay rate expected for an initially excited atom. (b) The total probability
  associated with the atom and the photon field.}
  \label{fig:ex2_infmode}
\end{figure}

We lastly consider the $p = 40$ case---this is again sufficient to achieve
convergence to high accuracy for the first five non-zero $a_n$---for $\xi_0 = 1$. The results are
given in Figure \ref{fig:ex2_infmode}, and are similar to those shown in
Figure \ref{fig:a_infmode}. After the initial excitation, the behavior of
the various Hermite modes comprising the atom amplitude is nearly
identical to that shown in Figure \ref{fig:a_infmode}a for the first
example. The behavior of the atomic and photonic contributions to the total
probability are similar, except in this case the photonic probability
contains contributions both from the incoming wavepacket and from the
field induced by the decaying atom.

\begin{figure}[t]
  \centering
  \includegraphics[width=\linewidth]{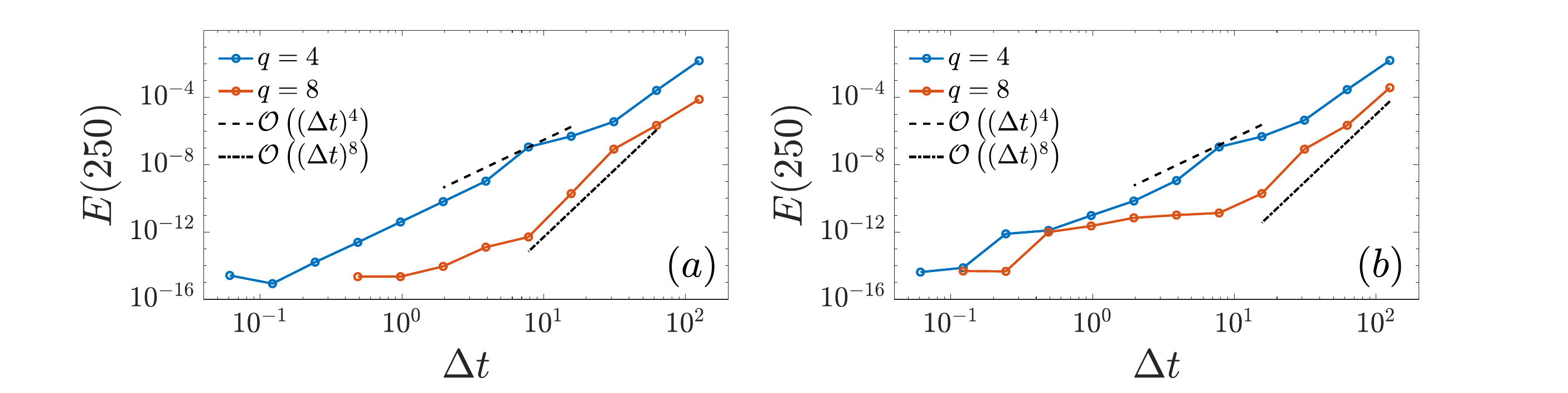}
  \caption{Error $E(t)$ for the second example with
  $\xi_0 = 1$ for (a) $p = 1$ and (b) $p = 40$, using fourth and eighth-order time-stepping.}
  \label{fig:ex2_err}
\end{figure}

We again verify the order of accuracy of the fourth and eighth-order
time-stepping algorithms by
measuring the error $E(t)$, defined by \eqref{eq:err}, for $t = 250$, with
$g = 0.2$ and $\xi_0 = 1$. Results for the $p = 1$ case are given in
Figure \ref{fig:ex2_err}a, and for the $p = 40$ case in Figure
\ref{fig:ex2_err}b.

\section{Conclusion} \label{sec:conclusion}

We have presented an efficient numerical method to solve
\eqref{eq:system} by reformulating it as an integro-differential
equation. This avoids the challenges associated with the nonlocality of
the differential operator, and the unboundedness of the domain. We
address the resulting Volterra-type memory dependence, for the case of a
Gaussian atomic density, by projecting the
solution history onto a collection of exponentials, which can be propagated
by a simple recurrence.

In our numerical experiments, when the spatial extent of the atom
amplitude is represented by a single degree of freedom, we recover the
expected Wigner-Weisskopf decay behavior for a one-atom system. When
multiple degrees of freedom are included, we observe more complicated
collective dynamics. Our numerical method serves as a useful starting
point to examine more complicated systems and related models in quantum
optics. In particular, in a forthcoming publication, we will generalize the method to 
systems of distinct two-level atoms coupled to a photon field.

\appendix

\section{Estimate of $j_n^{(2)}(t)$} \label{sec:jn2}

In this Appendix, we prove the estimate used to neglect $j_n^{(2)}(t)$
in Section \ref{sec:soe}. We have
\begin{align*}
  \frac{\Gamma\paren{\frac{n+1}{2}}}{2} \abs{j_n^{(2)}(t)} &= \abs{\int_0^\infty (\eta-ia)^n e^{-
(\eta-ia)^2 - i (\eta-ia) t} \, d\eta} \\
  &\leq e^{a^2-at} \int_0^\infty \paren{\eta^2+a^2}^{n/2} e^{-\eta^2} \,
  d\eta\\
  &= e^{2a^2-at} \int_a^\infty \frac{x^{n+1}}{\sqrt{x^2-a^2}}
  e^{-x^2} \, dx,
\end{align*}
where in the last line we have made the change of variables $x^2 =
\eta^2 + a^2$. We split the integral into two pieces:
\[\int_a^\infty \frac{x^{n+1}}{\sqrt{x^2-a^2}}
  e^{-x^2} \, dx = \paren{\int_a^{\sqrt{2}a} + \int_{\sqrt{2}a}^\infty} \frac{x^{n+1}}{\sqrt{x^2-a^2}}
  e^{-x^2} \, dx = I_1 + I_2.\]
For the first integral, we have
\[I_1 \leq 2^{(n+1)/2} a^{n+1} e^{-a^2} \int_a^{\sqrt{2} a}
\paren{x^2-a^2}^{-1/2} \, dx.\] 
The integral can be computed by the substitution $x = a \sec \theta$,
and is equal to $\log\paren{1+\sqrt{2}}$. We also use the estimate $a^{n+1} e^{-a^2}
\leq \sqrt{\frac{n+1}{2}}$ to obtain
\[I_1 \leq 2^{n/2} \sqrt{n+1} \log\paren{1+\sqrt{2}}.\]
For the second integral, we use that $\frac{x}{\sqrt{x^2-a^2}} \leq
\sqrt{2}$ when $x \geq \sqrt{2} a$ to obtain
\[I_2 \leq \sqrt{2} \int_{\sqrt{2}a}^\infty x^n e^{-x^2} \, dx \leq
\sqrt{2} \int_0^\infty x^n e^{-x^2} \, dx =
\frac{\Gamma\paren{\frac{n+1}{2}}}{\sqrt{2}}.\]
Combining these results gives the desired result,
\[\abs{j_n^{(2)}(t)} \leq e^{2a^2-at} \paren{\log\paren{1+\sqrt{2}}
\frac{2^{\frac{n}{2}+1} \sqrt{n+1}}{\Gamma\paren{\frac{n+1}{2}}} +
\sqrt{2}} \leq 14 e^{2a^2-at}.\]
In the last inequality, we have used that
$\frac{2^{\frac{n}{2}} \sqrt{n+1}}{\Gamma\paren{\frac{n+1}{2}}}$ reaches
its maximum of approximately $6.9$ at $n = 5$.

\section*{Acknowledgments}
John Schotland was supported in part by the NSF grant DMS-1912821 and
the AFOSR grant FA9550-19-1-0320. The Flatiron Institute is a division
of the Simons Foundation. 
  
\bibliographystyle{ieeetr}
\bibliography{qo1}

\end{document}